\documentclass[12pt,centertags,oneside]{amsart}

\usepackage{tikz}
\usepackage[mathscr]{eucal}
\usepackage{asymptote}
\usepackage{caption}
\usetikzlibrary{arrows}
\usepackage{import}
\usepackage{graphicx,graphics,array,booktabs}

\usepackage{amsmath,amstext,amsthm,amscd,typearea,hyperref}
\usepackage{amssymb}
\usepackage{a4wide}
\usepackage[mathscr]{eucal}
\usepackage{mathrsfs}
\usepackage{typearea}
\usepackage{charter}
\usepackage{pdfsync}
\usepackage{multirow}
\usepackage{array}
\usepackage{amsmath, amssymb, latexsym, amscd, amsthm,amsfonts,amstext}
\usepackage[mathscr]{eucal}

\usepackage{amsmath,amstext,amsthm,amscd,typearea,stmaryrd,mathabx}
\usepackage{amssymb}
\usepackage{a4wide}
\usepackage[mathscr]{eucal}
\usepackage{mathrsfs}
\usepackage{typearea}
\usepackage{charter}
\usepackage{pdfsync}
\usepackage[a4paper,width=16.2cm,top=3cm,bottom=3cm,lines=50]{geometry}

\numberwithin{equation}{section}

\newtheorem{theorem}{Theorem}[section]
\newtheorem{definition}[theorem]{Definition}
\newtheorem{proposition}[theorem]{Proposition}
\newtheorem{corollary}[theorem]{Corollary}
\newtheorem{lemma}[theorem]{Lemma}
\newtheorem{remark}[theorem]{Remark}

\newcommand{\cali}[1]{\mathscr{#1}}
\newcommand{\Tan}{\mathop{\mathrm{Tan}}\nolimits}

\newcommand{\comp}{{\mathop{\mathrm{comp}}\nolimits}}
\usepackage{graphicx} 

\newcommand{\AB}{{\rm AB}}

\newcommand{\Leb}{{\rm Leb}}

\newcommand{\Tube}{{\rm Tube}}
\newcommand{\volume}{{\rm vol}}

\newcommand{\supp}{{\rm supp}}

\newcommand{\ver}{{\rm ver}}
\newcommand{\hor}{{\rm hor}}

\renewcommand{\top}{{top}}
\newcommand{\ddc}{{dd^c}}

\newcommand{\ddcv}{{dd^c_\ver}}
\newcommand{\dc}{{d^c}}
\newcommand{\dbar}{{\overline\partial}}
\newcommand{\ddbar}{{\partial\overline\partial}}

\newcommand{\CL}{{\rm CL}}
\newcommand{\GL}{{\rm GL}}

\newcommand{\tot}{{\rm tot}}
\newcommand{\ind}{{\bf 1}}

\newcommand{\id}{{\rm id}}

\newcommand{\reg}{{\rm reg}}

\newcommand{\bfr}{{\rm \mathbf{r}}}

\newcommand{\bfi}{{\rm \mathbf{i}}}

\newcommand{\bfj}{{\rm \mathbf{j}}}
\newcommand{\bfU}{{\rm \mathbf{U}}}
\newcommand{\bfW}{{\rm \mathbf{W}}}

\newcommand{\upm}{{\mathrm{ \overline{m}}}}
\newcommand{\lowm}{{\mathrm{ \underline{m}}}}

\newcommand{\Cc}{\cali{C}}

\newcommand{\Ec}{\cali{E}}

\newcommand{\Kc}{\cali{K}}

\renewcommand{\Mc}{\cali{M}}

\newcommand{\Uc}{\cali{U}}
\newcommand{\Vc}{\cali{V}}

\newcommand{\Ic}{\cali{I}}

\newcommand{\FS}{{\rm FS}}

\newcommand{\Cf}{\mathfrak{C}}

\newcommand{\C}{\mathbb{C}}
\newcommand{\D}{\mathbb{D}}
\newcommand{\E}{\mathbb{E}}

\newcommand{\N}{\mathbb{N}}

\newcommand{\R}{\mathbb{R}}
\newcommand{\T}{\mathbb{T}}
\newcommand{\B}{\mathbb{B}}

\newcommand{\U}{\mathbb{U}}

\renewcommand{\P}{\mathbb{P}}

\title[Uniqueness  of tangent  currents]{Uniqueness  of tangent  currents for positive closed currents}

\author{Vi{\^e}t-Anh Nguy{\^e}n}
\address{Universit\'e de Lille, 
Laboratoire de math\'ematiques Paul Painlev\'e, 
CNRS U.M.R. 8524,  
59655 Villeneuve d'Ascq Cedex, 
France. }

\address{and Vietnam Institute for Advanced Study in Mathematics (VIASM),  157 Chua Lang Street, Hanoi, Vietnam.
}

\email{Viet-Anh.Nguyen@univ-lille.fr, {\tt   https://pro.univ-lille.fr/viet-anh-nguyen/}}

\author{Tuyen  Trung Truong}
\address{Department of Mathematics, University of Oslo, Blindern 0851 Oslo, Norway.
}
\email{ tuyentt@math.uio.no}

\date{February 10, 2025}

\begin{document}


\begin{abstract}
Let $X$ be a complex  manifold $X$ of dimension $k,$ 
and let $V\subset X$ be   a  K\"ahler submanifold  of dimension $l,$ and let   $B\subset V$ be  a piecewise $\Cc^2$-smooth  domain.  
  Let $T$ be a positive closed currents of bidegree $(p,p)$ in $X$  such that $ T$ satisfies a mild reasonable assumption  in a neighborhood of  $\partial B$ in $X$ and that
  the $j$-th average mean  $\nu_j( T, B,r)$ for every $j$  with $\max(0,l-p)\leq j\leq \min(l,k-p) $   converges  sufficiently  fast   to the  $j$-th generalized  Lelong number  $\nu_j(T, B )$ as $r$ tends to $0$ 
so that  $r^{-1}(\nu_j(T, B,r)-\nu_j( T,B))$  is  locally integrable near $r=0.$ 
Then  we  show that $T$  admits  a unique tangent current along $B.$  A  local version where we replace the condition of $T$ near $B$ by the conditions on a finite cover of $B$ by   piecewise $\Cc^2$-smooth  domains in $V$ is  also  given. 
When $T$ is a  current of integration over a complex analytic set,
we show that  $\nu_j( T, B,r)-\nu_j(T, B)=O(r^\rho)$  for some $\rho>0,$ and  hence this condition is  satisfied.
Our result may be viewed as a natural  generalization of Blel-Demailly-Mouzali's criterion from the case $l=0$ to the case   $l>0.$ 
The  result  has applications in the intersection theory of positive closed currents. 
\end{abstract}

\maketitle

\medskip\medskip


 \medskip



\section{Introduction} \label{S:Intro}

 Let $X$ be a  complex   manifold of dimension $k.$ Let   $d,$ $\dc$ denote the  real  differential operators on $X$  defined by
$d:=\partial+\overline\partial,$  $\dc:= {1\over 2\pi i}(\partial -\overline\partial) $ so that 
$\ddc={i\over \pi} \partial\overline\partial.$   
Let $V\subset  X$ be   a submanifold  of dimension $l.$
Let $T$ be a positive closed currents of bidegree $(p,p)$   on    $X.$
The existence  of tangent  currents  to $T$ along $V$  as well as their properties  have  recently received a lot of attention.   
In this  work we  address  the problem of finding  sufficient  conditions  ensuring  that 
there is  exactly one tangent  current. We begin  by recalling  the following classical case.


\subsection{Tangent  currents  for  $l=0$}

In this case  $V$ is reduced  to a single point $x\in X.$
Let $X$ be an  open  neighborhood of $x$ in $\C^k.$

Choose a local holomorphic coordinate system $z$ near $x$ such that $x=0$ in these coordinates. 
The Lelong number $\nu(T,x)$ of $T$ at $x$ is the limit of the normalized
mass of $\|T \|$  on the ball $\B(0, r)$ of center $0$ and radius $r$ when $r$ tends to $0.$ More
precisely, we have
\begin{equation}\label{e:Lelong-number-point}
\nu(T , x) := \lim\limits_{r\to 0}\nu(T,x,r),\quad\text{where}\quad \nu(T,x,r):=
{\sigma_T(\B(0,r))\over 
(2\pi)^{k-p} r^{ 2k-2p}}. 
\end{equation}
Here, $\sigma_T:={1\over (k-p)!}\, T\wedge ({i\over 2}\ddbar \|z\|^2)^{k-p}$ is  the trace measure of $T.$
Note that $(2\pi)^{k-p} r^{2k-2p}$  is the mass on $\B(0, r)$ of the $(p , p  )$-current of
integration on a linear subspace of dimension $k - p$ through $0.$  When $T$ is a positive closed current, Lelong establishes  in \cite{Lelong} (see also \cite{Lelong68})
that 
the  {\it  average mean} $\nu(T,x,r)$  is a  non-negative-valued  increasing function in the  radius $r.$
So  the limit \eqref{e:Lelong-number-point} always exists.


Let $A_\lambda : \C^k \to  \C^k$ be defined by
$A_\lambda (y) := \lambda y$ for $y\in\C^k$  with $\lambda \in \C^* .$  When $\lambda$ goes to infinity, the domain of definition
of the current $T_\lambda := (A_\lambda)_* (T)$ converges to $\C^k .$ This family of currents
is relatively compact,  and any limit current $T_\infty$ for $\lambda \to\infty,$ is  
called a {\it tangent current} to $T.$ 
A tangent  current  is defined on the  whole $\C^k,$ and it is 
conic  in the  sense  that  it is  invariant under $(A_\lambda)_* .$  The existence of   tangent currents  was initiated by many authors,  in particular 
by Harvey's exposition \cite{Harvey}.
Given a  tangent  current $T_\infty$  to $T,$  we can
extend it to $\P^k$ with zero mass on the hyperplane at infinity. Thus, there is a
positive closed current $\T_\infty$ on $\P^{k-1}$  such that $T_\infty=\pi^*_\infty
(\T_\infty).$ Here we identify the hyperplane at infinity with $\P^{k-1},$
and we denote by $\pi_\infty: \P^k \setminus \{0\} \to \P^{k-1}$
the canonical central projection. The class of $\T_\infty$ (resp. of $T_\infty$) in the de Rham
cohomology of $\P^{k-1}$ (resp. of $\P^k$) is equal to $\nu(T , x)$ times the class of a
linear subspace. So these cohomology classes do not depend on the choice of
$T_\infty.$ In  general, the tangent current $T_\infty$ is not unique, see  Kiselman   \cite{Kiselman}. 
The  following  result of   (see also Blel \cite{Blel89} for an earlier seminal result)
gives an effective condition  guaranteeing the uniqueness of the tangent currents. 
This is also the main motivation  for our work.

\begin{theorem}\label{T:BDM} {\rm (Blel-Demailly-Mouzali  \cite{BlelDemaillyMouzali})}
 Let $T$ be a  positive closed $(p,p)$ current on $X$  and $x$  a point in $X.$  Suppose that  one
 of  the  following  conditions  is  satisfied for  some $r_0>0.$
 \begin{equation}
 (i)\quad  \int_0^{r_0} {\sqrt{\nu(T,x,r) - \nu(T,x,r/2)}\over r}dr <\infty;\qquad  (ii)\qquad
  \int_0^{r_0} {\nu(T,x,r) - \nu(T,x)\over r}dr <\infty.
 \end{equation}
 Then $T$ admits a unique  tangent current.
\end{theorem}

\subsection{Tangent  currents     for $l > 0$ and   Dinh-Sibony theory,  and   Alessandrini-Bassanelli's  Lelong  number}

Next,  we deal with  the  situation where    the single point $x$ in the  previous  subsection  is  replaced  by   a submanifold $V\subset X$ of positive dimension $l$  ($1\leq l< k$).
Only recently,   Dinh and  Sibony \cite{DinhSibony18} have  developed  a  satisfactory  theory of tangent  currents and  density currents   for  positive  closed currents  in  this
context.

  Let $\E$ be the normal
vector bundle to $V$ in $X$ and $\pi:\ \E\to V$ be the  canonical projection.
Consider a point $x\in V.$ If $\Tan_x(X)$  and $\Tan_x(V)$ denote, respectively, the tangent spaces of $X$ and of $V$ at $x,$
the  fiber $\E_x$ of $\E$ over $x$ is canonically  identified  with the quotient space $\Tan_x(X)/\Tan_x(V).$

For $\lambda \in \C^\ast ,$ let $A_\lambda :\ \E \to  \E$ be the multiplication by $\lambda$ in fibers of $\E,$
that is, 
\begin{equation}\label{e:A_lambda} A_\lambda(y):=\lambda y\qquad\text{for}\qquad y\in \E.
\end{equation}
A  current $T$ on $\E$ is  said to be {\it $V$-conic} if $T$ is  invariant under  the action of  $A_\lambda,$ that is,
$(A_\lambda)_*T=T$ for all $\lambda\in\C^*.$
We identify $V$ with the zero section of $\E.$  

Let $\pi_0:\ \overline \E:=\P(\E\oplus\C)\to V$ be its canonical compactification. 

We  expect as  in  Harvey's exposition \cite{Harvey} that every  tangent  current  $T_\infty$ lives on $\E.$
However,  a basic difficulty arises.  When $V$ has
positive dimension, in general, no neighbourhood of $V$ in $X$ is biholomorphic
to a neighbourhood of $V$ in $\E.$

To  overcome  this  difficulty,  Dinh and  Sibony propose  a  softer  notion:  {\it the  admissible maps.}  More  precisely, let $\tau$  be a diffeomorphism between a neighbourhood of $V$ in $X$ and a
neighbourhood of $V$ in $\E$ whose restriction to $V$ is identity. We assume that
$\tau$ is admissible in the sense that the endomorphism of $\E$ induced by the
differential of $\tau$  when restricted to $V$ is the identity map from $\E$ to $\E.$

 Fix $0\leq p\leq  k$ and set
\begin{equation}\label{e:m}\upm:= \min(l,k-p)\qquad\text{and}\qquad
  \lowm:=\max(0,l-p).
  \end{equation}
Here is  the main result of  Dinh and  Sibony.
\begin{theorem}\label{T:Dinh-Sibony-first}{\rm  (Dinh-Sibony \cite[Theorems 1.1, 4.6 and Definition 4.8]{DinhSibony18})}  Let $X,$ $V,$ $\E,$  $\overline \E,$ $A_\lambda$ and $\tau$ be as above. 
Let $T$ be  a positive closed $(p,p)$-current on $X.$
Assume in addition  that $X$ is  K\"ahler and  $\supp(T)\cap V$ is  compact. Then:
\begin{enumerate} \item The family
of currents $T_\lambda:= (A_\lambda)_* \tau_* (T )$ is relatively compact and any limit current, for
$\lambda\to\infty,$ is a positive closed $(p, p)$-current on $\E$ whose trivial extension is a
positive closed $(p, p)$-current on $\overline \E.$   Such a limit current $S$ is  called a {\rm tangent current to $T$ along $V.$} 
\item If $S$ is a tangent current to $T$ along $V$, then it is
$V$-conic, i.e., invariant under $(A_\lambda)_* ,$  and its de Rham cohomology class $\{S\}$ in the cohomology group with compact support
$H^{2p}_\comp (\E, \C)$ does not depend on the choice of $\tau$ and $S.$
 
\end{enumerate}
\end{theorem}

To prove  their theorem,  Dinh and  Sibony develop a cohomological calculus on every positive closed current $T$ such that   $\supp(T)\cap V$ is  compact.
It is  worth noting that later on,  Vu \cite{Vu21} weakens the assumption of K\"ahlerian  on $X.$

 Dinh-Sibony's results  provide an effective tool  to measure the size of the intersection for arbitrary currents even when there is an excess of the intersection dimension. The
excess phenomenon cannot be handled  using the classical intersection theory for positive closed
currents.
 The reader may consult  Fulton's book  \cite{Fulton} for the algebraic  setting.

Dinh-Sibony  theory  has   many applications in Complex Analysis  and  Algebraic  Geometry, especially in  Complex  Dynamics and the theory of foliations.  
 See   \cite{DinhNguyenTruong15, DinhNguyenTruong17,DinhSibony16,Kaufmann,KaufmannVu,Vu21}  etc.
 
 On the other hand,  Alessandrini and Bassanelli introduce in  \cite{AlessandriniBassanelli96} a remarkable  notion of Lelong number  of   positive plurisubharmonic currents. Recall  that  a   real $(p,p)$-current  $T$ is  said to be  plurisubharmonic if the current  $\ddc T$ is positive. 
 So all positive closed currents are positive plurisubharmonic.
 In the context of    Alessandrini and Bassanelli,     
$(X,V)$  is a  very special pair of manifolds (affine manifolds), however  they  allow a domain $B\Subset V$ and formulate their Lelong number along $B.$
This means that $\supp(T)\cap B$ may be noncompact in $B.$

\begin{theorem}\label{T:AB-1}{\rm  (Alessandrini-Bassanelli \cite[Theorem I and Definition 2.2]{AlessandriniBassanelli96})}
Consider $X=\C^k$ and $V$ is a linear complex subspace of dimension $l\geq 0. $   
  We use the coordinates $(z,w)\in\C^{k-l}\times \C^l$ so  that  $V=\{z=0\}.$ Let $0\leq p<k-l$ and  
  let $T$  be    a positive plurisubharmonic $(p,p)$-current on an open  neighborhood $\Omega$ of $0$ in $\C^k.$ 
Then, for every open ball $B$ in $V,$ $B\Subset \Omega,$ the  following limit  exists and is  finite
$$
 \nu_\AB(T,B):=\lim_{r\to 0+} {1\over r^{2(k-l-p)}}\int_{\Tube(B,r)}  T(z,w)\wedge (\ddc \|z\|^2)^{k-l-p}\wedge (\ddc \|w\|^2)^{l},
 $$
 where the  tube  $\Tube(B,r)$ of radius $r$ over $B$ is  given by 
 \begin{equation} \label{e:Tube-AB}
\Tube(B,r):=\left\lbrace (z,w)\in \C^{k-l}\times\C^{l}:\ \|z\|<r,\ w\in B\right\rbrace .
\end{equation}
 $\nu_\AB(T,B)$ is  called    the {\rm  Alessandrini-Bassanelli's  Lelong  number of $T$ along  $B$}.
\end{theorem}

The  important viewpoint of Alessandrini-Bassanelli is  that when $V$ is of positive dimension,   tubular neighborhoods $\Tube(B,r)$ of $B$ and a  mixed form $  (\ddc \|z\|^2)^{k-l-p}\wedge (\ddc \|w\|^2)^{l}$ should replace
the usual balls $\B(x,r)$ around a single point $x$ with  the usual form  $  (\ddc \|z\|^2)^{k-p}.$
When $V$ is  a single point $\{x\}$ and $B=\{x\},$   Alessandrini-Bassanelli's  Lelong  number $\nu_\AB(T,x)$  coincides with the classical Lelong number
$\nu(T,x).$ 

\subsection{The generalized Lelong numbers}
The first named  author introduces in \cite{Nguyen21} a concept of  the Lelong numbers which generalizes  Alessandrini-Bassanelli's  Lelong  number 
in  a general setting of pairs of manifolds $(V,X).$  Instead of only one  Lelong number as in Theorem \ref{T:AB-1}, he defines a finite  set of cardinality $\leq  \dim V+1=l+1$ generalized Lelong numbers.   
 We recall here  the formulation of these characteristic  numerical quantities.
 
Let $X$ be a complex  manifold of dimension $k,$
 $V\subset  X$    a submanifold  of dimension $l,$ and $B\subset V$ a  relatively compact piecewise  $\Cc^2$-smooth open subset.
Let  $V_0$ be  a relatively compact  open subset of $V$ such that  $B\Subset V_0.$ Let $\omega$ be a Hermitian form on $V.$
Let $\tau:\ U\to\tau(U)$ be  an admissible  map along $B$
 on an open neighborhood $U$ of $\overline B$ in $X,$ see Definition \ref{D:admissible-maps} below.

  Denote by $\pi:\ \E\to V$ the canonical projection.
Consider a Hermitian metric  $h=\|\cdot\|$  on the  vector bundle  $\E_{\pi^{-1}(V_0)}$  and    let   $\varphi:\ \E_{\pi^{-1}(V_0)}\to \R^+$ be the function defined by  
\begin{equation}\label{e:varphi-spec}
 \varphi(y):=\|y\|^2\qquad \text{for}\qquad  y\in \pi^{-1}(V_0)\subset \E.
\end{equation}
Consider also the following  closed  $(1,1)$-forms on  $ \pi^{-1}(V_0)\subset \E $
\begin{equation}\label{e:alpha-beta-spec}
 \alpha:=\ddc\log\varphi\quad\text{and}\quad \beta:= \ddc\varphi.
\end{equation} 
So, for every $x\in V_0\subset  X$  the  metric $\| \cdot\|$  on the fiber $\E_x\simeq \C^{k-l}$ is  an Euclidean metric (in a suitable basis). In particular, we have 
\begin{equation}\label{e:varphi_bis-spec}  \varphi(\lambda y)=|\lambda|^2\varphi(y)\qquad\text{for}\qquad  y\in \pi^{-1}(V_0)\subset \E,\qquad\lambda\in\C. 
\end{equation}
For $r>0$ 
consider the  following {\it tube with base $B$ and radius $r$}
\begin{equation}
\label{e:tubular-nbh-0}
\Tube(B,r):=\left\lbrace y\in \E:\    \pi(y)\in B\quad\text{and}\quad  \|y\|<r  \right\rbrace.
\end{equation}
So   this is  a  natural generalization of Euclidean  tubes   considered by  Alessandrini-Bassanelli  in \eqref{e:Tube-AB}.
For for  all $0\leq s<r<\infty,$  define also the  corona tube
\begin{equation}\label{e:tubular-corona-0}\Tube(B,s,r):=\left\lbrace y\in \E:\   \pi(y)\in B\quad\text{and}\quad   s<\|y\|<r \right\rbrace.
\end{equation}
Since  $V_0\Subset V,$ there  is  a constant $c>0$ large enough such that 
$c\pi^*\omega+\beta$
is  positive on $\pi^{-1}(V_0).$  Moreover,   the latter  form defines
a K\"ahler metric there  if  $\omega$ is  K\"ahler on $V_0.$

 Let $\bfr>0$ be small enough such that $\Tube(B,\bfr)\subset  \tau(U),$  see \eqref{e:tubular-nbh-0}.
 Fix $0\leq p\leq k.$
 Let $T$ be a real  current of degree $2p$ and  of  order $0$    on    $U.$ 
 For $0\leq j\leq \upm$ and  $0<r\leq \bfr,$   consider 
\begin{equation}\label{e:Lelong-numbers}
 \nu_j(T,B,r,\omega,\tau,h):=  {1\over r^{2(k-p-j)}}\int_{\Tube(B,r)} (\tau_*T)\wedge \pi^*(\omega^j) \wedge \beta^{k-p-j}.  
\end{equation} 
When $j=\upm$  we  also denote  $\nu_\upm(T,B,r,\tau)$    by $\nu_\top(T,B,r,\tau).$

Let  $0\leq j\leq \upm.$ For  $0<s<r\leq \bfr,$   consider 
\begin{equation}\label{e:Lelong-corona-numbers}
 \kappa_j(T,B,s,r,\omega, \tau,h):=   \int_{\Tube(B,s,r)} (\tau_*T)\wedge \pi^*(\omega^j) \wedge \alpha^{k-p-j}. 
\end{equation} 
Let $0<r\leq\bfr.$   Consider
\begin{equation}\label{e:Lelong-log-bullet-numbers}
 \kappa^\bullet_j(T,B,r,\omega, \tau,h):= \limsup\limits_{s\to0+}   \kappa_j(T,B,s,r,\tau)
\end{equation} 
provided that the RHS side makes sense.
We  also consider
\begin{equation}\label{e:Lelong-log-numbers}
 \kappa_j(T,B,r,\omega,\tau,h):= \int_{\Tube(B,r)} (\tau_*T)\wedge \pi^*(\omega^j) \wedge \alpha^{k-p-j},  
\end{equation}
provided that the RHS side makes sense.

 \subsection{Classes    $\CL^{p;m,m'}$ ($0\leq m'\leq m$)  of positive closed currents  and a  tangent theorem}\label{SS:Classes}
 
 Let $X$ be a complex manifold of dimension $k.$
Let $V\subset  X$ be   a  K\"ahler submanifold  of dimension $l,$ and  $\omega$ a K\"ahler form on $V,$ 
and $B\subset V$
a  relatively compact piecewise  $\Cc^2$-smooth open subset.  Recall from \cite[Subsection 1.5]{Nguyen21}
the following   classes of positive closed currents  relevant  to our study.

\begin{definition}\rm \label{D:Class}
Let $m,m'\in\N$ with $0\leq m'\leq m.$
Let $T$ be a positive closed current of bidegree $(p,p)$ on $X.$
We say that {\it $T$ is  approximable along $B$ with $\Cc^m$-smooth positive  closed forms with $\Cc^{m'}$-control on boundary,}
if, there are an open neighborhood $U$ of $\overline B$ in $X,$   and an open neighborhood $W$ of $\partial  B$ in $X,$  and  a sequence of positive closed  forms $T_n\in\Cc^m(U)$     such that 
 \begin{enumerate}
\item the masses  $\|T_n\|$ on $U$ are uniformly bounded;

\item $T_n$ converge weakly  to $T$ on $U$ as $n$ tends to infinity;

\item     The $\Cc^{m'}$-norms of $T_n$'s on $W$ are uniformly bounded, that is,  
$$\sup_{n\geq 1} \|T_n\|_{\Cc^{m'}(W)} <\infty.$$
 \end{enumerate}
 $(T_n)_{n=0}^\infty$ is called an {\it approximating    sequence of $T$.}
 Let  $ \CL^{p;m,m'}(B)$ denote the class of all positive  closed currents on $X$ satisfying the above  properties.
 \end{definition}

 The relevance of the K\"ahlerian assumption on $X$ is illustrated by the  following result.
 
 \begin{proposition}{\rm  (See \cite[Theorem 1.19]{Nguyen21} which is  in turn inspired by \cite{DinhSibony04})}
   Suppose that $X$ is   K\"ahler  and let $m,m'\in\N$ with $0\leq m'\leq m.$ Then, for  every   positive closed current $T$ of bidegree $(p,p)$ on $X$  which is  of class $\Cc^{m'}$ in a neighborhood of $\partial B$ in $X$
 is of the form $T=T^+-T^-$ where $T^\pm\in \CL^{p;m,m'}(B).$
  \end{proposition}

\begin{remark}\label{R:approx-current-int} \rm In particular,  every current of integration $[S]$ on a complex analytic subvariety $S$ of pure codimension $p$ in $X$  such that  $S\cap B\Subset B$  belongs to  $\CL^{p;m,m'}(B)$ for every $m,m'\in\N$ with $0\leq m'\leq m.$
 \end{remark}

 The following  result  gives the existence of both tangent currents and   the generalized Lelong numbers with positive closed  currents. Here,    the notion of  strongly  admissible maps will be
 recalled   
  in  Definition \ref{D:Strongly-admissible-maps}.
 \begin{theorem}\label{T:Nguyen} {\rm  (Nguy\^en \cite[Tangent Theorem II (Theorem 1.11)]{Nguyen21})}
 Let $X$ be a complex  manifold of dimension $k.$
Let $V\subset  X$ be   a K\"ahler submanifold  of dimension $l,$  and $\omega$ a K\"ahler form on $V,$ 
and $B\subset V$
a  relatively compact piecewise  $\Cc^2$-smooth open subset.
   Let $T$ be  a   positive closed  $(p,p)$-current  on a neighborhood of $\overline B$ in $X$ such that       $T =T^+-T^-$ on a  neighborhood of $\overline B$ in $X,$
 where $T^\pm\in \CL^{p;1,1} (B).$
  Then  the following  assertions  hold:
   \begin{enumerate}
   \item For every $\lowm\leq j\leq \upm,$  the following limit  exists and is  finite
  $$    \nu_j(T,B,\omega,h):=\lim\limits_{r\to 0+}\nu_j(T,B,r,\omega,\tau ,h)                              $$
   for  all strongly  admissible maps  $\tau$     for $B$ and for all  Hermitian  metrics  $h$ on $\E.$   
   $    \nu_j(T,B,\omega,h)$ is  called {\rm  the $j$-th generalized  Lelong number of $T$  along $B$ associated  to $\omega$ and $h.$}
   
\item The  real numbers  $\nu_j(T,B,\omega,h)$ are intrinsic, that is, they are
  independent of  the choice  of   $\tau.$ 
  
  \item  The following  equality holds
  $$\lim\limits_{r\to 0+}\kappa_j(T,B,r,\tau,\omega,h)=\nu_j(T,B,\omega,h)$$ 
  for all $\lowm\leq j\leq \upm$ with $j>l-p,$ and  for all strongly  admissible maps $\tau$  for $B$ and for all Hermitian  metrics $h$ on $\E.$
  \item   $\nu_{\upm}(T,B,\omega,h)$ is  nonnegative. Moreover,  it is  totally intrinsic, i.e.  it is independent of  the choice  of both    $\tau$ and $h.$  So we will denote it simply by $\nu_{\upm}(T,B,\omega).$
  Moreover, it has  a  geometric meaning in the sense of Siu and Alessandrini-Bassanelli \cite{Siu,AlessandriniBassanelli96}.

   \item 
  There  exist   tangent currents  to  $T$  along $B$   and all  tangent  currents  $T_\infty$  are $V$-conic  positive  closed   on $\pi^{-1}(B)\subset \E.$
  
  \end{enumerate}
 \end{theorem}

The main idea   in \cite{Nguyen21} for the  proof of Theorem \ref{T:Nguyen}   is a systematic use of some  new  Lelong--Jensen type formulas for holomorphic  vector bundles, see Theorem \ref{T:Lelong-Jensen-closed} below  (which is \cite[Theorem 4.5]{Nguyen21}) for a particular  form of these formulas.
 This allows  for  treating   currents  $T$ where  $\supp(T)\cap V$ is  noncompact in $V.$

\begin{remark}\rm
 Since the forms $\alpha$ and $\beta$ are in general not positive and the strongly admissible map $\tau$ is  in general not holomorphic,  the $\nu_j(T,B,\omega,h)$'s for $\lowm\leq j<\upm$  are in general  not nonnegative. 
 
 The  work in \cite[Theorem 2.17]{Nguyen24a} provides a sufficient condition  ensuring that $    \nu_j(T,B,\omega,h)$ is totally intrinsic for $\lowm\leq j\leq \upm.$  
\end{remark}

\subsection{Statement of the  main results}

Recall that $l:=\dim V.$
The theory of tangent currents  for $l=0$ has    Theorem \ref{T:Dinh-Sibony-first} as its analogue for $l>0.$
The classical  Lelong number   for $l=0$  has   Theorems \ref{T:AB-1} and \ref{T:Nguyen} as its analogues for $l>0.$
Therefore, we may  expect  that Theorem  \ref{T:BDM} has also an  analogue for $l>0.$

Our  following  first main  result gives an effective  sufficient condition ensuring  the uniqueness of the tangent currents.
It may be regarded as  a counterpart of  Theorem  \ref{T:BDM} for the case where $l >0.$
\begin{theorem}\label{T:main_1} Let $X$ be a complex  manifold of dimension $k.$
Let $V\subset  X$ be   a K\"ahler submanifold  of dimension $l,$ $\omega$ a K\"ahler form on $V,$ 
and $B\subset V$
a  relatively compact piecewise  $\Cc^2$-smooth open subset. 
Let $\tau$ be a  strongly admissible map for $B.$
  Let $h$  be a Hermitian metric  on $\E|_{\overline B}.$ 
Let $T$ be a positive closed current of bidegree $(p,p)$ on $X$  
 such that      $T =T^+-T^-$ on a  neighborhood of $\overline B$ in $X,$
 where $T^\pm\in \CL^{p;1,1} (B),$ and that  
        one of the  following two conditions is  fulfilled:
  \begin{itemize}
  \item [(a)] $$\int_0^{r_0} {\sqrt{|\nu_j(T, B,r,\tau, h)-\nu_j(T, B,r/2,\tau,h)|}\over  r}dr <\infty\qquad \text{for all}\qquad \lowm\leq j\leq \upm;$$
  \item [(b)]  $$ \int_0^{r_0} {\big |\nu_j(T,B,r,\tau,h)-\nu_j(T, B,h)\big|\over  r}dr <\infty\qquad \text{for all}\qquad \lowm\leq j\leq \upm.$$
  \end{itemize}
for some $r_0>0$  small  enough.
Then $T$ admits a  unique  tangent  current along $B.$  
\end{theorem}

Our second main result  provides a relevant ``local'' sufficient  condition for    the uniqueness of the tangent currents.
It is very convenient in practice, see  \cite{Nguyen24a} for  applications in the intersection theory of positive closed currents.

\begin{theorem}\label{T:main_2} Let $X$ be a complex manifold of dimension $k.$
Let $V\subset  X$ be   a  K\"ahler submanifold  of dimension $l,$ and  $\omega$ a K\"ahler form on $V,$ 
and $B\subset V$
a  relatively compact piecewise  $\Cc^2$-smooth open subset admitting a  strongly admissible map. 
Suppose in addition   that   $\overline B$ can be covered  by a finite collection of    $\Cc^2$-smooth  domains $(B_i)_{i\in I}$ with $B_i\Subset V$  such that $B_i\cap B$ is  a  $\Cc^2$-smooth open subset in $V.$  Let $\tau_i$ be a  strongly admissible  map for $B_i,$ and $h_i$   a Hermitian metric  on $\E|_{\overline B_i}$ for $i\in I.$
Let $T$ be a positive closed currents of bidegree $(p,p)$  
 such that     
     $T =T^+-T^-$ on a  neighborhood of $\overline B$ in $X,$
 where $T^\pm\in \CL^{p;1,1} (B),$ and that 
    for every $i\in I,$   one of the  following two conditions is  fulfilled:
  \begin{itemize}
  \item [(a-i)] $$\int_0^{r_0} {\sqrt{|\kappa_j(T,B_i\cap B,r/2, r,\tau_i,h_i)|}\over  r}dr <\infty\qquad\text{ for all}\qquad \lowm\leq j\leq \upm;$$
  \item [(b-i)]  $$ \int_0^{r_0} {|\kappa^\bullet_j(T,B_i\cap B,r,\tau_i,h_i)|\over  r}dr <\infty \qquad\text{ for all}\qquad \lowm\leq j\leq \upm    ;$$
  \end{itemize}
for some $r_0>0$  small  enough.
Then $T$ admits a  unique  tangent  current along $B$.  
\end{theorem}

We give here  the  following consequence   which captures the  essential points of the above  two main results in the  special but important context where the ambient manifold $X$ is K\"ahler. 
\begin{corollary}\label{C:main}
 Let $X,V$ be as above   and  suppose that $X$ is K\"ahler and that $(V,\omega)$ is K\"ahler  and  that $B\subset V$ is 
a  relatively compact piecewise  $\Cc^2$-smooth open subset.  Let $T$ be a  positive closed current of bidegree $(p,p)$ on $X$ such that there is a neighborhood $W$ of $\partial B$ in $X$ such that the restriction of $T$ on $W$ is    a $\Cc^1$-smooth form. Then $T$  admits a unique  tangent current  along $B$ if at least one of the  following item $\bullet$  is  fulfilled:
\begin{itemize}
 \item either condition (a) or condition (b) in Theorem \ref{T:main_1} is  satisfied;
 
 \item  for every $i\in I,$   either  condition (a-i) or condition (b-i) in Theorem \ref{T:main_2} is  satisfied.
\end{itemize}
\end{corollary}

 \begin{remark}
  \rm In \cite{Vu21}  under the assumption that $\supp T\cap V\Subset V,$ Vu  proves the  existence of tangent  currents which are $V$-conic (see assertion (5) of Theorem  \ref{T:Nguyen}), using  a Hermitian
  metric  form $\hat\omega$ on $X$  such that  $\ddc \hat\omega^j=0$ on $V$ for $1\leq j\leq \upm-1.$ 
  Instead of  the K\"ahlerian assumption  on $(V,\omega),$  we may assume  the following weaker assumption $\ddc \omega^j=0$ on $V$ for $1\leq j\leq \upm-1.$  Note that this  assumption is  weaker than that of Vu \cite{Vu21} since  $\omega$ is  only defined on $V$ whereas $\hat\omega$ is defined on $X.$    However, we do require  the existence of a strongly admissible map whereas Vu only needs  admissible ones. 
  Moreover, we  assume that $T\in\CL^{p;1,1}(B).$ Then  the conclusions of Theorems \ref{T:main_1} and \ref{T:main_2} still hold. However, we do not enter into details for such a  result.
 \end{remark}

Our last main result  studies 
the particular example of the current of integration over a complex analytic set $S\subset X.$
\begin{theorem}\label{T:main_3} Let $X$ be a K\"ahler  manifold of dimension $k.$
Let $V\subset  X$ be   a  submanifold  of dimension $l,$ $\omega$ a K\"ahler form on $V,$ 
and $B\subset V$
a  relatively compact piecewise  $\Cc^2$-smooth open subset. 
Let $\tau$ be a   admissible map for $B.$
  Let $h$  be a Hermitian metric  on $\E|_{\overline B}.$ Let $S$ be  a complex  analytic set of pure codimension $p$ in $X$ such that $S\cap B\Subset B,$  
  and let   $[S]$ be its current of integration. 
Then there are constants $\rho, c>0$  such that 
  $$ |\nu_j([S],B,r,\tau,h)-\nu_j([S], B,h)|\leq cr^\rho \qquad \text{for all}\qquad \lowm\leq j\leq \upm,\qquad 0<r\leq \bfr.$$
   In particular,  by Theorem \ref{T:main_1} (b) $[S]$ admits a  unique  tangent  current along $B.$  
\end{theorem}

Although  the uniqueness of the  tangent current   
in this  situation is already well-known at least when $S\cap V\Subset V$ (see e.g \cite[Proposition 2.6]{Vu21}), Theorem \ref{T:main_3}
may  provide  a  quantified speed of  the 
 convergence of $(A_\lambda)_*(\tau_*T)$ as $\lambda$ tends to infinity.

The article is  organized as follows. 

In Section \ref{S:preparation}  we set up the background of our study.  More precisely, we recall some basic facts of the theory of positive currents, the definition of 
 (strongly) admissible maps and a Lelong-Jensen formula for closed currents on holomorphic vector bundles established in \cite{Nguyen21}.

In Section \ref{S:Global-estimates}  we first  establish   a key  global inequality  
(Proposition \ref{P:Kc_j,q})   which is  a  consequence of   a basic mass estimate (Theorem \ref{T:basic-estimate}) obtained in \cite{Nguyen21}.
Next, using  the technique of admissible estimates as developed in \cite{Nguyen21} (see also \cite{Nguyen24a}) and the above  key global inequality,
we  reduce the conditions  in Theorems  \ref{T:main_1} and \ref{T:main_2} to an appropriate  growth behavior of  some global mass indicators in Propositions \ref{P:main-1} and \ref{P:main-2}.

Next, Section \ref{S:Local-estimates} infer from the  previous global estimates  corresponding local     growth behavior of appropriate local mass indicators.

Having  good local estimates at hands, we  apply, in Section \ref{S:Proofs}, the limiting argument of Blel-Demailly-Mouzali \cite{BlelDemaillyMouzali} in order to get the 
desired convergence of $(A_\lambda)_*(\tau_*T)$ as $\lambda$ tends to infinity.

Using the  construction of the cylindrical cones in \cite{Nguyen24b}, Section \ref{S:Cone}    gives the proof of Theorem \ref{T:main_3}.
\smallskip

\noindent
{\bf Acknowledgments. }   The   authors are supported by the Labex CEMPI (ANR-11-LABX-0007-01).
The  first named author is supported by
  the project QuaSiDy (ANR-21-CE40-0016).
The paper was partially prepared 
during the visits of the first named  author at the Vietnam  Institute for Advanced Study in Mathematics (VIASM) and 
at the University of Oslo,  and during  the  visit of the  second  named author  at  the University of Lille. They would like to express their gratitude to these organizations for hospitality and  for  financial support.

\section{Preparatory results}\label{S:preparation}

\subsection{Positive  currents}\label{SS:currents}

 Let $T$ be  a current of bidegree $(p,p)$ on an open set $U\subset\C^k.$  Write 
\begin{equation}\label{e:expression_currents}T=i^{p^2}\sum T_{I,J} dx_I\wedge d\overline x_J\quad\text{with}\quad T_{I,J}\quad \text{a distribution on}\quad U,
\end{equation}
the sum being taken all  over all  multi-indices $I,J$ with $|I|=|J|=p.$
Here, for a multi-index $I=(i_1,\ldots,i_p)$  with $1\leq i_1\leq \ldots\leq i_p\leq k,$  $|I|$ denotes the length $p$ of $I,$   $dx_I$ denotes $dx_{i_1}\wedge \ldots  dx_{i_p}$  and
 $d\bar x_I$ denotes $d\bar x_{i_1}\wedge \ldots  d\bar x_{i_p}.$

\begin{proposition}\label{P:Demailly}{\rm (see  e.g.  \cite[Proposition 1.14]{Demailly}) }
 Let  $T=i^{p^2}\sum T_{I,J} dx_I\wedge d\overline x_J$ of bidegree $(p,p)$ be  a  positive  current on an open  set in $\C^k.$  Then its coefficients $T_{I,J}$ are complex measures
 and satisfy  $\overline T_{I,J}=T_{J,I}$ for all  multi-indices $|I|=|J|=p.$  Moreover,  $T_{I,I}\geq 0,$ and the absolute values $|T_{I,J}|$
 of the measure $T_{I,J}$ satisfy the  inequality
 $$
 \lambda_I\lambda_J|T_{I,J}|\leq 2^{k-p}\sum_M  \lambda_M^2 T_{M,M},  \qquad  I\cap J\subset M\subset I\cup J,
 $$
 where $\lambda_j\geq 0$   are  arbitrary   coefficients and $\lambda_I=\prod_{j\in I}   \lambda_j.$
\end{proposition}

 Let $R$ be  a current with measure coefficients (or equivalently,  of order $0$) on an open set $\Omega$  in a complex  manifold $X$ of dimension $k.$
 Let $W$ be a relatively compact open subset of $\Omega$ and $\Phi$ a smooth test form  on $\Omega,$ we will  write
 \begin{equation}\label{e;cut-off}
 \int_W  R\wedge \Phi:=  \langle  R,\ind_W\Phi\rangle,
 \end{equation}
 where $\ind_W$ is  the  characteristic  function of $W.$ Let $(R_n)_{n=1}^\infty$ be a  sequence of positive  currents on $\Omega$ such that
 $\lim_{n\to\infty}  R_n=R$  weakly on  $\Omega,$ then  we  see that
 \begin{equation}\label{e:continuity-cut-off}
 \lim_{n\to\infty}\int_W  R_n\wedge \Phi=\int_W  R\wedge \Phi
 \end{equation}
 for every smooth test form $\Phi$  on $\Omega$ and  every relatively compact open subset $W\subset \Omega$ with $\|R\|(\partial W)=0.$
 Here, $\partial W$ is  the topological boundary of $W$ and $\|R\|$ is  the mass-measure of $R.$
 Consequently, if $K$  is a compact subset of $\Omega$ and $(W_i)_{i\in I}$ is a family of open subsets of $\Omega$ such that
 $K\subset  W_i$ for all $i\in I$ and $\partial W_i\cap \partial W_{j}=\varnothing$ for $i\not=j,$ then  we have
 \begin{equation} \label{e:except-countable}
   \lim_{n\to\infty}\int_{W_i}  R_n\wedge \Phi=\int_{W_i}  R\wedge \Phi
 \end{equation}
for every smooth test form $\Phi$  on $\Omega$ and  every $i\in I$  except  for a countable subset of $I.$

 The following  notions of   quasi-positive currents and their  weak convergence are needed.
 \begin{definition}{\rm  (\cite[Definition 2.3]{Nguyen21})}  \label{D:quasi-positivity}\rm
  We say that   a    current $R$ defined on $\Omega$  is {\it  quasi-positive} if, for every $x\in \Omega,$ there are an open neighborhood $\Omega_x$ of $x$  in $\Omega$ and  a $\Cc^1$-diffeomorphism $\tau_x$
 of $\Omega_x$ such that $\tau_x^*R$ is a  positive current.
 
 We say that   a  sequence of  currents $(R_n)_{n=1}^\infty$ {\it converge  in the sense of quasi-positive currents  on $\Omega$ to a current $R$}
  if   for every $x\in \Omega,$ there are an open neighborhood $\Omega_x$ of $x$  in $\Omega$ and  a $\Cc^1$-diffeomorphism $\tau_x$
 of $\Omega_x$  and two sequences of  positive  currents $(T^\pm_n)_{n=1}^\infty$  on $\Omega_x$   such that all currents $\tau_x^*(R_n-R)= T^+_n-T^-_n$   and that both sequences $T_n^\pm$ converge  weakly to  a common positive current  $T$ on $\Omega_x.$  
 \end{definition}
 
 Here is  a  simple  property  of these notions.
 \begin{lemma}{\rm  (\cite[Lemma 2.4]{Nguyen21})}\label{L:quasi-positivity}
  If  a  sequence of  currents $(R_n)_{n=1}^\infty$  converge  in the sense of quasi-positive currents  on $\Omega$ to a current $R,$ then 
  both  \eqref{e:continuity-cut-off} and \eqref{e:except-countable}  hold.
 \end{lemma}

\subsection{(Strongly) admissible maps}\label{SS:Admissible}

The  following notion,  introduced  by  Dinh-Sibony \cite[Definitions 2.15 and  2.18]{DinhSibony18}, plays a vital role in their tangent theory for positive closed currents.  
\begin{definition}\label{D:admissible-maps}   \rm   
 Let $B$ be a relatively compact nonempty open subset of $V.$
An admissible map   along  $B$  is   a $\Cc^1$-smooth diffeomorphism $\tau$  from an open neighborhood  $U$ of    $\overline{B}$ in $X$
onto an open neighborhood of $B\subset V$ in $\E$ (where $V$ is  identified with the zero section $0_\E$) which is identity on an open neighborhood of  $\overline B\subset V$ such that the  endomorphism   on $\E$ induced  by the restriction
of the differential $d\tau$ to ${\overline B}$ is identity.  

In local coordinates, we can describe an admissible
map $\tau$ as follows: for every  point $x\in V\cap U,$  for every local
chart $y=(z,w)$  on a neighborhood $W$ of   $x$ in $U$  with  $V\cap W=\{z=0\}$, 
 we have
\begin{equation}\label{e:admissible-maps}
\tau (y) = \big( z + O(\|z\|^2), w+ O(\|z\|\big) ,
\end{equation}
and
\begin{equation}
d\tau  (y) = \big( dz + \widetilde O(\|z\|^2), dw + \widetilde O(\|z\| ) \big),
\end{equation}
as $z \to  0$ where for every positive integer $m,$  $\widetilde O(\|z \|^m )$ denotes the sum of $1$-forms with
$O(\|z\|^m )$-coefficients and a linear combination of $dz,$ $d\bar z$  with $O(\|z\|^{m-1} )$-coefficients.
\end{definition}
Observe that in  \cite{DinhSibony18}   Dinh-Sibony  use the terminology {\it almost-admissible}  for those
maps satisfying  Definition \ref{D:admissible-maps}.
In general, $\tau$ is not holomorphic. When $U$ is a small enough local chart,
we can choose a holomorphic  admissible map by using suitable holomorphic coordinates
on $U .$   For the  global situation, the following  result gives a positive answer.
\begin{theorem}{\rm (\cite[Lemma 4.2]{DinhSibony18})}
For every  compact  subset $V_0\subset V,$ there always exists  an admissible map $\tau$ defined on a small enough tubular neighborhood  $U$ of $V_0$ in $X.$
\end{theorem}
     
     The first named  author,  inspired
by Dinh-Sibony \cite[Proposition 3.8]{DinhSibony18b}, has  introduced  in \cite[Definition 2.7]{Nguyen21} the following notion, which allows him 
to  develop   a quantitative theory of  tangent  and density currents for  positive plurisubharmonic currents. 
\begin{definition}
 \label{D:Strongly-admissible-maps}\rm  Let $B$ be a relatively compact nonempty open subset of $V.$
   A {\it strongly admissible} map  along  $B$ is   a $\Cc^2$-smooth diffeomorphism $\tau$  from an open neighborhood $U$ of $\overline B$  in $X$
onto an open neighborhood of $V \cap U$ in $\E$  such that for every  point $x\in V\cap U,$  for every local
chart $y=(z,w)$  on a neighborhood $W$ of   $x$ in $U$  with  $V\cap W=\{z=0\}$,
 we have   
 \begin{eqnarray*}
 \tau_{j}(z,w)&= &  z_j+\sum_{p,q=1}^{k-l}  a_{pq}(w) z_pz_{q}  +O(\|z\|^3)\quad\text{for}\quad 1\leq j\leq  k-l,\\
 \tau_{j}(z,w)&= &  w_{j-(k-l)} +\sum_{p=1}^{k-l} b_p(w)z_p 
 +O(\|z\|^2)\quad\text{for}\quad k-l< j\leq  k.
 \end{eqnarray*}
 Here,  we  write $\tau(y)=(\tau_{1}(y),\ldots,\tau_{k-l}(y),\tau_{k-l+1}(y),  \ldots \tau_{k}(y))\in\C^k,$  and $a_{pq},$  $b_p,$ 
 are $\Cc^2$-smooth functions depending only on $w.$
 In other words, if we write $\tau(z,w)=(z',w')\in\C^{k-l}\times \C^l,$ then
 \begin{eqnarray*}
 z'&= &  z+  z A z^T  +O(\|z\|^3),\\
 w'&= &  w + Bz 
 +O(\|z\|^2),
 \end{eqnarray*}
 where  $A$ is a $(k-l)\times(k-l)$-matrix and $B$ is  a $l\times(k-l)$-matrix whose entries are $\Cc^2$-smooth functions in $w,$ $z^T$ is the transpose of $z,$
\end{definition}

It is  worth noting that  a  strongly  admissible map is necessarily  admissible in the sense of Definition \ref{D:admissible-maps}.
 On the  other hand,  holomorphic  admissible maps  are always strongly   admissible, see \cite[Subsection 2.2]{Nguyen21}.

 The relevance of the K\"ahlerian assumption on $X$ is illustrated by the  following result.
 
 \begin{proposition}{\rm  (See \cite[Theorem 1.19]{Nguyen21} which is  in turn inspired by \cite{DinhSibony18b})}
   Suppose that $X$ is   K\"ahler.
Then there always exists a strongly admissible  map $\tau$ along $B.$
  \end{proposition}
 
 \subsection{Lelong-Jensen formula for closed currents}\label{SS:Lelong-Jensen}
 
 According to  \cite[Section 4]{Nguyen21}, for   all $r_0\leq r<s\leq \bfr,$  define
\begin{equation}\label{e:tubular-corona}\Tube(B,r,s):=\left\lbrace y\in \E:\   \pi(y)\in B\quad\text{and}\quad   r^2<\varphi(y)<s^2  \right\rbrace.
\end{equation}
Note that the boundary   $\partial\Tube(B,r)$ can be decomposed as the disjoint union  of the {\it  vertical boundary}   $\partial_\ver\Tube(B,r)$ and
the  {\it horizontal  boundary}    $\partial_\hor\Tube(B,r)$, where
\begin{eqnarray*}
 \partial_\ver\Tube(B,r)&:=&  \left\lbrace y\in \E:\   \pi(y)\in \partial B\quad\text{and}\quad  \varphi(y)\leq r^2  \right\rbrace ,\\
 \partial_\hor\Tube(B,r)&:=&  \left\lbrace y\in \E:\   \pi(y)\in B\quad\text{and}\quad  \varphi(y)=r^2  \right\rbrace .
\end{eqnarray*}
Recall  from   \cite[Theorem 4.5]{Nguyen21}   the following Lelong-Jensen formula for  currents  which  are  approximable by smooth {\bf closed}  forms with  control on the boundary.

\begin{theorem}\label{T:Lelong-Jensen-closed}
Let $r\in(0,\bfr]$ and  let $S$ be a real closed current of dimension $2q$ on  a neighborhood of $\overline\Tube(B,r).$ 
%
%
Suppose that there is a sequence of $\Cc^1$-smooth closed forms of dimension $2q:$  $(S_n)_{n=1}^\infty$ defined on  a neighborhood of $\overline\Tube(B,r)$ 
such  that
$S_n$ converge to $S$  in the sense of quasi-positive currents   on a neighborhood of  $\overline\Tube(B,r)$ as $n$ tends to infinity (see Definition \ref{D:quasi-positivity}).

Then the following  two assertions  hold:
\begin{enumerate} \item  The  following two  sub-assertions hold:
For  all $r_1,r_2\in  [0,r]$ with  $r_1<r_2$  except for  a  countable set of   values, we have that  
\begin{equation}\label{e:Lelong-Jensen-closed}
 {1\over  r_2^{2q}} \int_{\Tube(B,r_2)} S\wedge \beta^q- {1\over  r_1^{2q}} \int_{\Tube(B,r_1)} S\wedge \beta^q
 =\lim\limits_{n\to\infty} \Vc(S_n,r_1,r_2)+   \int_{\Tube(B,r_1,r_2)} S\wedge \alpha^q.
\end{equation} 
Here the vertical boundary term  $\Vc(S,r_1,r_2)$ for  a  continuous form $S$  is  given by
\begin{equation}\label{e:vertical-boundary-term-closed}
\begin{split}
 \Vc(S,r_1,r_2)&:=
 {1\over r_2^{2q}} \int_{\partial_\ver \Tube(B,r_2)}\dc\varphi\wedge S\wedge \beta^{q-1}-{1\over r_1^{2q}} \int_{\partial_\ver \Tube(B,r_1)}\dc\varphi\wedge S\wedge \beta^{q-1}\\
  &-\int_{\partial_\ver \Tube(B,r_1,r_2)}\dc\log\varphi\wedge S\wedge \alpha^{q-1}.
  \end{split}
\end{equation}
\item If $S$ is a closed $\Cc^1$-smooth form, then identity \eqref{e:Lelong-Jensen-closed} (with $S_n:= S$ for $n \geq 1$) holds
for  all $r_1,r_2\in  (0,r]$ with  $r_1<r_2.$
\end{enumerate}
\end{theorem}


 \begin{remark}\label{R:comparison-and-slicing} \rm By Lemma \ref{L:quasi-positivity},  there is only a countable set of
$r_1 ,$ $ r_2$ for which the mass of $S$ on the boundary of $ \Tube(B, r_1 , r_2)$  is strictly positive. Therefore, 
 the  countable  set of values  appearing in Theorem   \ref{T:Lelong-Jensen-closed}   is obtained.
 
We  investigate  Theorem \ref{T:AB-1} by using  the  canonical metric $\|z\|^2$ for $\E=\C^{k-l},$ and the canonical metric
$\omega(w):=\|w\|^2$  for $w\in V=\C^l.$  We obtain that  $\beta:=\ddc \|z\|^2,$  $\alpha:=\ddc \log\|z\|^2.$
In this  context and  when   $S=T\wedge\pi^*\omega^l,$ where $T$ is a positive closed  current,   Theorem \ref{T:Lelong-Jensen-closed} yields Lelong-Jensen  type formulas in \cite[Proposition 2.5]{AlessandriniBassanelli96} and  the vertical boundary term  $\Vc(S,r_1,r_2)$  given  in  \eqref{e:vertical-boundary-term-closed}  vanishes. These formulas are essential   in order to   obtain  Theorem \ref{T:AB-1}.

As already observed  in \cite[Section 4]{Nguyen21}, when $S=T\wedge \pi^*\omega^j$ where  $T$ is a positive closed  current and  $\omega$ is  a K\"ahler form  on $V,$      the vertical boundary term  $\Vc(S,r_1,r_2)$  given  in  \eqref{e:vertical-boundary-term-closed}  vanishes  if 
$j=l,$ but if $j<l,$ it does not vanish in general. That is  why  one cannot use  the slicing method as in \cite {AlessandriniBassanelli96}
in order to obtain   Theorem \ref{T:Nguyen}.

When $l=0,$    $B=V=\{\text{a  single  point}\ x\in X\}$ and     Theorem   \ref{T:Lelong-Jensen-closed} yields  the classical Lelong-Jensen formula for positive closed  currents. 
 \end{remark}

\section{Technical tools}\label{S:Tools}

\subsection{Forms $\alpha_\ver$ and $\beta_\ver$}
 
 Since  the  transition  functions  of the holomorphic vector bundle $\E$ are holomorphic, 
 the vertical operators $\partial_\ver,$ $\dbar_v$ 
which are the  restrictions of the usual operators $\partial$ and $\dbar$ on fibers of $\E$  are well-defined. More precisely, for a smooth
form $\Phi$ on an open set $\Omega$ in $\E,$ we can define  as in \cite[Section 5.1]{Nguyen21}:
\begin{equation}\label{e:d-and-dbar-ver}
 \partial_\ver\Phi(y):=\partial |_{\E_{\pi(y)}}\Phi (y) \qquad\text{and}\qquad \dbar_\ver\Phi(y):=\dbar|_{\E_{\pi(y)}}\Phi(y)\qquad\text{for}\qquad y\in \Omega.
\end{equation}
So  the vertical operators  $d_\ver$ and  $\ddcv$ are  also well-defined  by the formulas \begin{equation}
\label{e:ddc-ver}
d_\ver \Phi:=\partial_\ver \Phi+\dbar_\ver \Phi\quad\text{and}\quad \ddcv\Phi:={i\over \pi}  \partial_\ver\dbar_\ver\Phi.                           
                          \end{equation}  
Consider for $ y\in \E,$
\begin{equation}\label{e:alpha-beta-ver}
 \alpha_\ver(y):=\ddcv \log\varphi (y)=\ddc|_{\E_{\pi(y)}}\log\varphi (y) \qquad\text{and}\qquad \beta_\ver(y):=\ddcv \varphi(y)=\ddc|_{\E_{\pi(y)}}\varphi(y),
\end{equation}
where  $\ddc|_{\E_{\pi(y)}}$  is restriction  of  the operator $\ddc$ on the fiber $\E_{\pi(y)}.$  
Note that  both
 $\alpha_\ver$  and $\beta_\ver$  are positive  $(1,1)$-forms on $\E.$
However, they are not necessarily closed.

\subsection{Analysis  in local coordinates}
  The following local model studied in   \cite[Section 5.1]{Nguyen21} is   useful. Consider  an open set $V_0\Subset  V$ and  let us  study  near a given point $y_0\in V_0.$
We use the coordinates $(z,w)\in\C^{k-l}\times \C^l$  around a neighborhood $U$ of $y_0$  such that $y_0=0$  in  these coordinates. 
We may assume that  $U$ has the form $U=U'\times U'',$ where $U'$ (resp. $U'')$ is an  open neighborhood of $0'$ in $\C^{k-l}$ (resp. of  $0''$ in $\C^l$),
and  $V=\{z=0\}\simeq U''.$ Moreover,  we  may assume  that $U''=(2\D)^l.$  
Consider  the trivial  vector bundle $\pi:\ \E \to  U''$ with  $\E\simeq  \C^{k-l}\times U''.$ 
There is a smooth function  $A:\  \D^l\to \GL(\C,k-l)$ such that
\begin{equation}\label{e:varphi-new-exp}
\varphi(z,w)= \| A(w)z\|^2\qquad\text{for}\qquad  z\in\C^{k-l},\ w\in \D^l,
\end{equation}
where $\varphi$ is   defined in \eqref{e:varphi-spec}.
It follows  from \eqref{e:alpha-beta-ver} and  \eqref{e:varphi-new-exp} that  
\begin{equation}\label{e:tilde-alpha-beta-local-exp}
\alpha_\ver(z,w)=  A(w)^* [\ddc \log{\|z\|^2}]\quad\text{and}\quad    \beta_\ver(z,w)=  A(w)^* [\ddc \|z\|^2]\quad\text{for}\quad  z\in\C^{k-l},\ w\in \D^l.
\end{equation}

\subsection{Strong Standing Hypothesis}\label{SS:Ex-Stand-Hyp}

In this subsection we  recall from  \cite[Section  7.1]{Nguyen21} a  standard   setting for further technical  developments.

Consider  a  strongly admissible map $\tau:\ \bfU\to\tau(\bfU)$  along $B,$  with $\bfU$ a neighborhood of $\overline B$ in $X.$
By shrinking $\bfU$ if necessary, we may   fix a finite collection $\Uc=(\bfU_\ell,\tau_\ell)_{1\leq \ell\leq \ell_0} ,$
 of holomorphic admissible maps for $\bfU.$  More precisely,  
there is a  finite cover of $\overline \bfU$ by open subsets $\bfU_\ell,$ $1\leq \ell\leq  \ell_0,$ of $X$
such that   there is   a holomorphic coordinate system on $\overline \bfU_\ell$ in $X$ and $\bfU_\ell$  is  biholomorphic  to  $\U_\ell:=\tau_\ell(\bfU_\ell)\subset \E$
by a  holomorphic admissible map $\tau_\ell.$ By  choosing $\bfr>0$ small enough, we may assume  without loss of generality that $\overline\Tube(B,\bfr)\Subset \U:=\bigcup_{\ell=1}^{\ell_0} \U_\ell.$ 
Fix  a partition of unity  $(\theta_\ell)_{1\leq \ell\leq \ell_0}$ subordinate   to the open cover  $(\bfU_\ell\cap V)_{1\leq \ell\leq \ell_0}$   of $\overline{\bfU\cap V}$  in $V$  such that $\sum_{1\leq \ell\leq \ell_0}  \theta_\ell=1$ on an open neighborhood of $\overline {\bfU\cap V} \subset V.$ We  may  assume  without loss of generality  that there are open  subsets
$\widetilde V_\ell\subset V$  for $1\leq \ell\leq \ell_0$ such that
\begin{equation}\label{e:supp} \supp(\theta_l)\subset  \widetilde V_\ell\Subset \bfU_\ell\cap V \quad\text{and}\quad 
\tau( \widetilde V_\ell)\Subset  \U_\ell\quad\text{and}\quad \pi^{-1}(\supp(\theta_\ell))\cap \U\subset \U_\ell.
\end{equation}
 For $1\leq \ell\leq \ell_0$ set 
 \begin{equation}\label{e:tilde-tau_ell}\tilde \tau_\ell:=\tau\circ\tau_\ell^{-1}.
 \end{equation}
  So  $\tilde \tau_\ell$ defines a map  from  $\U_\ell\subset \E$  onto  $\tau(\bfU_\ell)\subset \E.$

  
   
\subsection{Representative currents}\label{SS:Regularizations}

Recall that $\pi:\  \E\to V$ is  the canonical projection.
\begin{definition}\label{D:T-hash}\rm

Let  $ T$ be a current  defined on $\bfU.$ Consider the current  $T^{\#}$ defined  on $\U$ by the   following formula:
\begin{equation}\label{e:T-hash}
T^{\#}:=\sum_{\ell=1}^{\ell_0}  (\pi^*\theta_\ell)\cdot (\tau_\ell)_*( T|_{\bfU_\ell}).
\end{equation}
By  \eqref{e:supp},  $T^{\#}$ is  well-defined. Moreover, if $T$ is  positive then so is  $T^{\#}.$

Let   $0\leq s<r\leq \bfr.$  Consider the currents  $T^{\#}_r$  and $T^{\#}_{s,r}$  defined  on $\U$ as  follows:
\begin{equation}\label{e:T-hash_r}
T^{\#}_r:=\sum_{\ell=1}^{\ell_0}  (\pi^*\theta_\ell)\cdot (\ind_{\Tube(B,r)}\circ \tilde \tau_\ell )\cdot(\tau_\ell)_*( T|_{\bfU_\ell})\quad\text{and}\quad
T^{\#}_{s,r}:=\sum_{\ell=1}^{\ell_0}  (\pi^*\theta_\ell)\cdot (\ind_{\Tube(B,s,r)}\circ \tilde \tau_\ell )\cdot(\tau_\ell)_*( T|_{\bfU_\ell}).
\end{equation}
Here, for $K\subset \E,$ $\ind_K$  is  the characteristic  function associated to  $K,$ that is, $$\ind_K(y):=\begin{cases}
                                       1,&\text{if $x\in K$};\\
                                       0,& \text{otherwise}.
                                      \end{cases}
$$    Note that if $T$ is  positive then so are  $T^{\#}_r$ and $T^{\#}_{s,r}.$ 
\end{definition}
Assume that $T$ is a positive current. Since   $T^{\#}_r$ and $T^{\#}_{s,r}$
are positive, they are easier to work with  than $\tau_*T$  (which is not positive in
general). Hence, some analogs of the generalized  Lelong numbers for $T^{\#}_r$ and $T^{\#}_{s,r}$
are always non-negative. On
the other hand, there are some close relations between these variants   and the
generalized Lelong numbers of $T .$  Details about these facts are  discussed in the remaining of this section.

\subsection{Positive forms $\hat\alpha$ and $\hat\beta$}
Recall from \cite[Lemma 5.1]{Nguyen21}  the construction of   positive  currents/forms   $\hat\alpha,$ $\hat\alpha'$  and $\hat\beta.$  
This positivity plays  a  crucial  role    in the sequel.

\begin{lemma}\label{L:hat-alpha-beta} 
\begin{enumerate}
\item  There is  a  constant $c_1>0$ large enough such that 
\begin{equation}\label{e:hat-beta}
  \hat\beta:=c_1\varphi\cdot  \pi^*\omega+\beta
 \end{equation}
is  positive  on $\pi^{-1}(V_0)$ and is strictly positive on $\pi^{-1}(V_0)\setminus V_0,$
and 
\begin{equation}\label{e:hat-alpha'}
\hat\alpha':= c_1 \pi^*\omega+ \alpha
\end{equation}
satisfies 
\begin{equation}\label{e:hat-alpha'-vs-alpha_ver}
\hat\alpha'\geq  c^{-1}_1\alpha_\ver.
\end{equation}
In particular,   $\hat\alpha'$ is positive  on $\pi^{-1}(V_0).$
\item 
  There are constants $c_2,c_3>0$ such that on $\Tube(V_0,\bfr)\setminus V_0,$  
\begin{equation}\label{e:hat-alpha}
\hat\alpha:=\hat\alpha'+c_2\beta= c_1 \pi^*\omega+ \alpha+c_2\beta
\end{equation}
is   strictly positive, and 
\begin{equation}\label{e:hat-alpha-vs-alpha_ver}
\hat\alpha\geq  c^{-1}_1\alpha_\ver,
\end{equation}
and
\begin{equation}\label{e:hat-alpha-vs-hat-beta}
\varphi\hat\alpha\leq  c_3\hat\beta.
\end{equation}
\item 
  There are constants $c_3>0$ such that  on $\Tube(V_0,\bfr),$
  \begin{equation}\label{e:hat-beta-vs-beta_ver}
\hat\beta\geq  c^{-1}_1\beta_\ver,
\end{equation}
  and on $\Tube(V_0,\bfr)\setminus V_0,$  
\begin{equation}\label{e:tilde-alpha-vs-hat-beta}
\varphi\alpha_\ver\leq c_3\hat\beta.
\end{equation}
\end{enumerate}

\end{lemma}

\subsection{Local and global mass indicators for positive currents}

Recall  from  \cite[Section 8.1]{Nguyen21} 
  the following  mass indicators   for a  positive  current $T$  of bidegree $(p,p)$ defined on $X.$
 For $0\leq j\leq k$ and $0\leq q\leq k-l$ and  $1\leq\ell\leq \ell_0,$ and  for $0<s<r\leq \bfr,$ 
 \begin{equation}\label{e:local-mass-indicators-bis}
 \begin{split}
 \Mc_j(T,r,\tau_\ell)&:={1\over  r^{2(k-p-j)}}\int (\ind_{\Tube(B,r)}\circ \tilde\tau_\ell)  (\pi^*\theta_\ell)\cdot (\tau_\ell)_*( T|_{\bfU_\ell})
 \wedge\pi^*\omega^j\wedge  \hat\beta^{k-p-j},\\
 \Kc_{j,q}(T,r,\tau_\ell)&:=\int(\ind_{\Tube(B,r)}\circ \tilde\tau_\ell)  (\pi^*\theta_\ell)\cdot (\tau_\ell)_*( T|_{\bfU_\ell})
 \wedge\pi^*\omega^j\wedge  \hat\beta^{k-p-q-j}\wedge \hat \alpha^q,\\
 \Kc_{j,q}(T,s,r,\tau_\ell)&:=\int (\ind_{\Tube(B,s,r)}\circ \tilde\tau_\ell) (\pi^*\theta_\ell)\cdot (\tau_\ell)_*( T|_{\bfU_\ell})
 \wedge\pi^*\omega^j\wedge  \hat\beta^{k-p-q-j}\wedge \hat \alpha^q.
 \end{split}
 \end{equation}
\begin{remark}\rm As was already observed in  \cite[Corollary 4.8]{Nguyen21} 
that $(\tau_\ell)_*( T|_{\bfU_\ell})
 \wedge\pi^*\omega^\upm$ is of full bidegree $(l,l)$ in $\{dw,d\bar w\}.$
Consequently,
by the bidegree reason, we  deduce  that  $\Mc_j(T,r,\tau_\ell),$ $\Kc_{j,q}(T,r,\tau)$ and $\Kc_{j,q}(T,s,r,\tau)$
are equal to $0$ provided that $j>\upm.$ 
\end{remark}
 
  We define   the following  global mass indicators.
 \begin{equation}\label{e:global-mass-indicators}
 \begin{split}
 \Mc_j(T,r)=\Mc_j(T,r,\Uc)&:=\sum_{\ell=1}^{\ell_0}\Mc_j(T,r,\tau_\ell) ,\\
 \Mc^\tot(T,r)=\Mc^\tot(T,r,\Uc)&:=\sum_{j=0}^\upm \Mc_j(T,r),\\
 \Kc_{j,q}(T,r)=    \Kc_{j,q}(T,r,\Uc) &:=\sum_{\ell=1}^{\ell_0}\Kc_{j,q}(T,r,\tau_\ell),\\
 \Kc_{j,q}(T,s,r)=    \Kc_{j,q}(T,s,r,\Uc) &:=\sum_{\ell=1}^{\ell_0}\Kc_{j,q}(T,s,r,\tau_\ell).
 \end{split}
 \end{equation}

 \begin{lemma} {\rm (\cite[Lemma 8.2]{Nguyen21}) }
 \label{L:Mc_j} The following  equalities hold:
  \begin{equation*}
 \begin{split}
 \Mc_j(T,r)&= {1\over r^{2(k-p-j)}}\int  T^\hash_r\wedge\pi^*\omega^j\wedge  \hat\beta^{k-p-j} ,\\
 \Kc_{j,q}(T,r) &=     \int  T^\hash_r\wedge
 \pi^*\omega^j\wedge  \hat\beta^{k-p-q-j}\wedge \hat\alpha^q ,\\
 \Kc_{j,q}(T,s,r) &=\int  T^\hash_{s,r}
 \wedge\pi^*\omega^j\wedge  \hat\beta^{k-p-q-j}\wedge \hat\alpha^q.
 \end{split}
 \end{equation*}
  \end{lemma}
\proof
It follows  from \eqref{e:global-mass-indicators}, \eqref{e:local-mass-indicators-bis},   and \eqref{e:T-hash_r}.
\endproof

\begin{lemma}  {\rm (\cite[Lemma 8.3]{Nguyen21}) }  \label{L:comparison-Mc} For every constant $\rho >1$ there is a constant $c>0$ such that
for $0\leq j\leq k$ and for $0<r<s<\rho r\leq \bfr$ and for every  positive current $T$ of bidegree $(p,p)$ on $\bfU,$    we have
$\Mc_j(T,r)\leq c\Mc_j(T,s).$ In particular, it holds that $\Mc^\tot(T,r)\leq c\Mc^\tot(T,s).$
\end{lemma}

According to   \cite[Section 8.4]{Nguyen21},
consider the following mass indicators, for $\lowm \leq j\leq \upm:$ 
 \begin{equation}\label{e:other-mass-indicator} \hat\nu_j(T,r):={1\over  r^{2(k-p-j)}}\int\limits_{\Tube(B,r)}
 \tau_* T\wedge (\beta+c_1r^2\pi^*\omega)^{k-p-j}\wedge  \pi^*\omega^j.
 \end{equation}
 We also consider the following new mass indicator, where   $T^\hash_r$ is given  in \eqref{e:T-hash_r}:
\begin{equation}\label{e:other-mass-indicators-bis}
 \Mc^\hash_j(T,r):={1\over r^{2(k-p-j)}}\int
 T^\hash_r \wedge (\beta+c_1r^2\pi^*\omega)^{k-p-j}\wedge  \pi^*\omega^j
 .
 \end{equation}

The only difference between $ \beta+c_1r^2\pi^*\omega$ and $\hat \beta$  is that  
function $\varphi$ in the latter expression is  replaced by $r^2$ in the former one. 

\begin{lemma} \label{L:Mc-hash-vs-nu}{\rm  (\cite[Lemma 8.10]{Nguyen21})} There is a  constant $c>0$ such that for every $\lowm\leq j\leq \upm$ and $0<r\leq\bfr:$
$$|\Mc^\hash_j(T,r)- \hat\nu_j(T,r)| 
\leq    cr\sum^{\upm}_{q=\lowm}  \Mc^\hash_q(T,r).$$
\end{lemma}

 \begin{lemma}\label{L:comparison-Mc-bis}{\rm  (\cite[Proposition  8.11]{Nguyen21})}
 For $\lowm \leq j\leq \upm,$  we have that $$ \lim_{r\to 0+} \Mc^\hash_j(T,r)=\lim_{r\to 0+}  \hat\nu_j(T,r)=\sum\limits_{q=0}^{k-p-j}  {k-p-j\choose q} c^q_1\nu_{j+q}(T,B,\tau).$$
 \end{lemma}

 \subsection{Technique of admissible  estimates}\label{SS:admissible-technique}
 
 We recall some  elements of the  technique of admissible  estimates which has been developed in \cite[Sections 7 and  8]{Nguyen21}
 (see also \cite[Section 5]{Nguyen24a} for a refinement). 

Let $\omega$ be  a K\"ahler metric  on $V.$
Let $\bfj=(j_1,j_2,j_3,j_4)$ with $j_1,j_3,j_4\in\N $ and $j_2\in{1\over 4}\N,$ and 
$k-p-j_1-j_3\geq 0.$ For $0<s<r\leq\bfr,$ and  for a real current $T$  on $\bfU,$  consider  
\begin{equation}\label{e:I_bfj}\begin{split}  I_{\bfj}(s,r)&:=\int_{\Tube(B,s,r)}\tau_*T\wedge \varphi^{j_2}(c_1-c_2\varphi)^{j_4}\hat\beta^{k-p-j_1-j_3}\wedge
(\pi^*\omega)^{j_3}\wedge \hat\alpha^{j_1},\\ 
I^\hash_{\bfj}(s,r)&:=\int_{\Tube(B,s,r)}T^\hash_{s,r}\wedge \varphi^{j_2}(c_1-c_2\varphi)^{j_4}\hat\beta^{k-p-j_1-j_3}\wedge
(\pi^*\omega)^{j_3}\wedge \hat\alpha^{j_1}.
\end{split}
\end{equation}
We define  $I_\bfj(r)$ and $I^\hash_\bfj(r)$ similarly  replacing the current $T^\hash_{s,r}$ (resp.  the domain of integration $\Tube(B,s,r)$) by
$T^\hash_r$  (resp. $\Tube(B,r)$).
\begin{remark}\label{R:I_bfj}
 \rm  Observe that $\Kc_{j,q}(T,r)=I^\hash_{(q,0,j,0)}(r)$ and $\Kc_{j,q}(T,s,r)=I^\hash_{(q,0,j,0)}(s,r).$
\end{remark}

Recall  from \cite[Lemma 8.6]{Nguyen21} the  following technical  key result.
\begin{lemma}\label{L:spec-wedge}  There are   constants $c,c_0$ independent of $T$    such that the   following inequality holds for  $0<s<r<\bfr:$
 \begin{equation*}
|I_\bfj(s,r)- I^\hash_\bfj(s,r)|^2 \leq c\big(\sum_{\bfj'} I^\hash_{\bfj'}(s-c_0s^2,r+c_0r^2)\big)\big ( \sum_{\bfj''} I^\hash_{\bfj''}(s-c_0s^2,r+c_0r^2) \big).  
\end{equation*}
Here, on the RHS:
\begin{itemize} \item[$\bullet$] the first sum  is taken over a finite number of multi-indices    $\bfj'=(j'_1,j'_2,j'_3,j'_4)$ as above  such that  $j'_1\leq  j_1$  and $j'_2\geq j_2;$ and either ($j'_3\leq j_3$) or ($j'_3>j_3$ and $j'_2\geq j_2+{1\over 2}$).
\item  the second sum   is taken over  a finite number of multi-indices $\bfj''=(j''_1,j''_2,j''_3,j''_4)$ as above   such that   either  ($j''_1< j_1$)
or ($j''_1=j_1$ and $j''_2\geq {1\over 4}+j_2$) or ($j''_1=j_1$ and $j''_3<j_3$).
\end{itemize}
\end{lemma}

\begin{remark}\label{R:j_4-not-necessary}\rm 
 The index $j_4$  is not indispendable in the sequel since $ c_1-c_2\varphi =O(1)$ on $\Tube(B,\bfr) .$ However, this index  is  necessary for the proof of   \cite [Theorem 8.7]{Nguyen21}, which is 
 Theorem
 \ref{T:basic-estimate} below. Moreover, it helps to clarify some arguments. 
\end{remark}

\section{Global  core estimates}\label{S:Global-estimates}
  

Fix  an open neighborhood $\bfW$ of $\partial B$ in $X$ with $\bfW\subset \bfU.$

\begin{definition}\label{D:sup}\rm
Fix an open neighborhood $\bfU$ of $\overline B$ and an open neighborhood $\bfW$ of $\partial B$ in $X$ with $\bfW\subset \bfU.$
Let $\widetilde\CL^{p;1,1}(\bfU,\bfW)$ be the  set of all $T\in \CL^{p;1,1}(\bfU,\bfW)$  whose  a sequence of approximating  forms $(T_n)_{n=1}^\infty\subset \Cc^1(\bfU)$
satisfies the following   condition:
 \begin{equation}\label{e:unit-CL-1,1} \|T_n\|_{\bfU}\leq  1\qquad\text{and}\qquad  \| T_n\|_{\Cc^1(\bfW)}\leq 1.\end{equation}
 \end{definition} 
 Recall from \cite[Theorem 8.7]{Nguyen21} the following basic mass estimate.
 \begin{theorem}\label{T:basic-estimate}
  There is a constant $c>0$ such that  $\Kc_{j,q}(T, \bfr)<c$  all $T\in \widetilde\CL^{p;1,1}(\bfU,\bfW).$
 \end{theorem}

\subsection{Estimates of the mass indicator $\Kc_{j,q}$  using the technique of admissible estimates}

The  following  global inequality
plays a key role in the sequel.

\begin{proposition}\label{P:Kc_j,q} 
Let $T\in \widetilde\CL^{p;1,1} (B).$
Let $0\leq p\leq k$ and $\lowm \leq j\leq \upm$ and  $0\leq q\leq k-p-j.$
Then  there is a constant $c>0$ independent of $T$ such that  $\Kc_{j,q}(T,0, r)\leq  cr^{2(k-p-q -j)}$ for all $0<r\leq\bfr.$
\end{proposition}
\proof
If $q=k-p-j,$ then   the  result  follows  from   Theorem \ref{T:basic-estimate}.
Therefore, in the remainder of the proof,  we assume  that   $q<k-p-j.$
By  \eqref{e:hat-beta} and  \eqref{e:hat-alpha-vs-hat-beta}, we infer  that
$$\hat\alpha\leq c_1\pi^*\omega+ c_2\varphi^{-1}\beta.$$
Using this,  it follows that
\begin{eqnarray*}
\Kc_{j,q}(T,0,r) &= &    \int_{\Tube(B,r)}  T^\hash_r\wedge
 \pi^*\omega^j\wedge  \hat\beta^{k-p-q-j}\wedge \hat\alpha^q\\
 &\leq &  \int_{\Tube(B,r)}  T^\hash_r\wedge
 \pi^*\omega^j\wedge  (\beta+c_1\varphi \pi^*\omega)^{k-p-q-j}\wedge (c_1\pi^*\omega+c_2 \varphi^{-1}\beta )^q\\
 &\leq&  \sum_{n=0}^\infty\int_{\Tube(B,2^{-n-1}r, 2^{-n}r)}  T^\hash_r\wedge
 \pi^*\omega^j\wedge  (\beta+c_1   \varphi  \pi^*\omega)^{k-p-q-j}\wedge  (c_1\pi^*\omega+c_2\varphi^{-1}\beta )^q.
 \end{eqnarray*}
 Since by  \eqref{e:varphi-spec}, $\varphi(y)=\|y\|^2,$  we infer that $(r')^2\leq\varphi(y)\leq r^2$ for  $y\in\Tube (B,r',r)$  with $0<r'<r\leq\bfr,$ and hence  the  expression in the last line is dominated  by a constant times
 \begin{eqnarray*}
  && \sum_{s=0}^q\sum_{n=0}^\infty  (2^{-n}r)^{-2(q-s)}  \int_{\Tube(B, 2^{-n}r)}  T^\hash_r\wedge
 \pi^*\omega^j\wedge  (\beta+c_1  (2^{-n}r)^{2}\pi^*\omega)^{k-p-q-j}\wedge  (\pi^*\omega)^s\wedge (\beta )^{q-s}\\
 &\lesssim&
 \sum_{s=0}^q\sum_{n=0}^\infty  (2^{-n}r)^{-2(q-s)}  \int_{\Tube(B, 2^{-n}r)}  T^\hash_r\wedge
 \pi^*\omega^{j+s}\wedge  (\beta+c_1 (2^{-n}r)^{2}\pi^*\omega)^{k-p-j-s}.
\end{eqnarray*}
By  formula \eqref{e:other-mass-indicators-bis}, the above expression is equal to
\begin{equation*}
  \sum_{s=0}^q\sum_{n=0}^\infty 
  (2^{-n}r)^{2(k-p-j-q)} \Mc^\hash_{j+s}(T,2^{-n}r).
 \end{equation*}
Since by Lemmas 
 \ref{L:Mc-hash-vs-nu}  and \ref{L:comparison-Mc-bis} and Theorem \ref{T:Nguyen}, $ \Mc^\hash_{j+s}(T,r)\leq c$ independent of $0<r\leq \bfr$ and $T,$ the above  sum is  bounded by a constant
 times    $\sum_{s=0}^q\sum_{n=0}^\infty 
  (2^{-n}r)^{2(k-p-j-q)},$ which gives us the desired conclusion since  $ k-p-j-q>0$ by our assumption. Hence
 the result follows.
\endproof
 Next, using the
technique of admissible estimates   which is recalled  in  Subsection \ref{SS:admissible-technique}, and the above key
global inequality  
(Proposition \ref{P:Kc_j,q}),  the conditions in Theorems \ref{T:main_1} and \ref{T:main_2} are reduced to an appropriate
growth behavior of some global mass indicators in Propositions \ref{P:main-1} and \ref{P:main-2} below.
\begin{lemma}\label{L:j_4=0} Let $T\in \CL^{p;1,1} (B).$
 For $\bfj=(j_1,j_2,j_3,j_4)$ and  $\bfj'=(j_1,j_2,j_3,0),$ there is a constant $c>0$ such that
  $I^\hash_\bfj(s,r)\leq c  I^\hash_{\bfj'}(s,r)$ for $0\leq s<r\leq \bfr.$
\end{lemma}
\proof
It follows from the observation that since  by \eqref{e:varphi-spec} $\varphi(y)=\|y\|^2,$ ${c_1\over 2}\leq c_1-c_2\varphi (y) <2c_1$ when $y\in \Tube(B,r)$ and $0<r\leq \bfr.$ 
\endproof

\begin{proposition}\label{P:main} Let $T\in \CL^{p;1,1} (B).$ Let $0\leq p\leq k$ and $\lowm \leq j\leq \upm$ and set $\bfj:=(k-p-j,0,j,0).$
Let $0\leq s<r\leq \bfr.$
Then   $|I_\bfj (s,r) - I^\hash_\bfj(s,r)|\leq  cr^{1/4}.$
\end{proposition}
\proof
Applying  Lemma 
\ref{L:spec-wedge} yields  a constant $c$ independent of $T$ and $s,r$   such that the   following inequality holds
 \begin{equation*}
|I_\bfj(s,r)- I^\hash_\bfj(s,r)|^2 \leq c\big(\sum_{\bfj'} I^\hash_{\bfj'}(s-c_0s^2,r+c_0r^2)\big)\big ( \sum_{\bfj''} I^\hash_{\bfj''}(s-c_0s^2,r+c_0r^2) \big).  
\end{equation*}
Here, on the RHS:
\begin{itemize} \item[$\bullet$] the first sum  is taken over a finite number of multi-indices    $\bfj'=(j'_1,j'_2,j'_3,j'_4)$ as above  such that  $j'_1\leq  j_1$  and $j'_2\geq j_2;$ and either ($j'_3\leq j_3$) or ($j'_3>j_3$ and $j'_2\geq j_2+{1\over 2}$).
\item  the second sum   is taken over  a finite number of multi-indices $\bfj''=(j''_1,j''_2,j''_3,j''_4)$ as above   such that   either  ($j''_1< j_1$)
or ($j''_1=j_1$ and $j''_2\geq {1\over 4}+j_2$) or ($j''_1=j_1$ and $j''_3<j_3$).
\end{itemize} 
By Lemma \ref{L:j_4=0} and the  first item $\bullet$, the  first sum $\sum_{\bfj'} I^\hash_{\bfj'}(s-c_0s^2,r+c_0r^2)$ is bounded by a constant times $\sum_j \Kc_{j,q}$ because $\varphi$ is  bounded on $\Tube(B,\bfr).$ 
By Proposition \ref{P:Kc_j,q}  this  sum is 
 bounded by a constant independently of $T.$
 
 By Lemma \ref{L:j_4=0} and the  second item $\bullet$, the  second  sum  is bounded by a constant times 
 $$\varphi^{1/4}\sum_j \Kc_{j,q}+  \sum_{j: j<k-p-q} \Kc_{j,q}$$ 
By Proposition \ref{P:Kc_j,q}  the second sum in the last line  is bounded by a constant times $r^2.$ 
In all, the    second sum  $\sum_{\bfj''} I^\hash_{\bfj''}(s-c_0s^2,r+c_0r^2)$ is bounded by a constant times $r^{1\over 2}.$
 This, combined with  the previous estimate on   $\sum_{\bfj'} I^\hash_{\bfj'}(s-c_0s^2,r+c_0r^2),$  implies the result.
\endproof

\begin{proposition}\label{P:main-bis}  There is a constant $c>0$  such that for every
   $T\in \widetilde\CL^{p;1,1} (\bfU,\bfW)$ and $0\leq s<r\leq \bfr,$ the following   assertions hold:
  \begin{enumerate}
   \item 
  $$    c^{-1}\big(\sum_{j=\lowm}^\upm \Kc_{j,k-p-j}(T, B,s,r) +r^{1\over 4} \big)\leq  \sum_{j=\lowm}^\upm |  \kappa_j(T, B,s, r,\tau,h)|+ r^{1\over 4}  \leq c\big(\sum_{j=\lowm}^\upm \Kc_{j,k-p-j}(T, B,s,r) +r^{1\over 4} \big).$$
\item For every $i\in I,$
\begin{equation*}\begin{split}   & c^{-1}\big(\sum_{j=\lowm}^\upm \Kc_{j,k-p-j}(T, B_i\cap B,s,r) +r^{1\over 4} \big)\leq  \sum_{j=\lowm}^\upm |  \kappa_j(T,B_i\cap  B,s, r,\tau,h)|+ r^{1\over 4} \\
&\leq c\big(\sum_{j=\lowm}^\upm \Kc_{j,k-p-j}(T, B_i\cap B,s,r) +r^{1\over 4}  \big).\end{split}\end{equation*}
  \end{enumerate}
 \end{proposition}
 \proof Following  the proof of \cite[Theorem 8.7]{Nguyen21},
consider
\begin{equation}\label{e:P-Lc-bullet_finite(0)}
 \Kc^{-}(T,s,r):=\sum\limits_{j+q<k-p}\Kc_{j,q}(T,s,r),
 \quad
 \Kc_{q}(T,s,r):=\sum\limits_{q'\leq q}\Kc_{j,q'}(T,s,r).
\end{equation}
We  define   $\Kc^{-}(T,r)$  and   $\Kc_{q}(T,r)$ similarly.
 Let $0\leq q_0\leq \upm$ and  set $j_0:= k-p-q_0.$
Recall from \eqref{e:hat-alpha} and \eqref{e:hat-beta} that
\begin{equation*} \alpha=\hat\alpha- c_1 \pi^*\omega-c_2\beta
=\hat\alpha -c_2\hat \beta +(c_2\varphi-c_1)\pi^*\omega
\quad\text{and}\quad  \beta=\hat\beta-c_1\varphi\cdot  \pi^*\omega.
 \end{equation*}
So we get  that  
\begin{eqnarray*} 
 \alpha^{q_0}
&=&(\hat\alpha -c_2\hat \beta +(c_2\varphi-c_1)\pi^*\omega)^{q_0} \\
&=&
\hat \alpha^{q_0}
+\sum_{j_1,j'_1}^q {q_0 \choose j_1}{q_0-j_1 \choose j'_1}(-c_2)^{j_1'}\cdot\,\hat\beta^{j'_1}\wedge 
((c_2\varphi-c_1)\pi^*\omega)^{q_0-j_1-j'_1}\wedge \hat\alpha^{j_1},
\end{eqnarray*}
where the last sum is  taken over all $(j_1,j'_1)$   such that   $0\leq j_1,j'_1\leq q_0$  and
$j_1+j'_1\leq  q_0$ and 
$j_1\not=q_0.$  Using this 
and the first equality of  \eqref{e:I_bfj}, we  have 
\begin{equation}\label{e:P-Lc-bullet_finite(5)}
 \begin{split}
\kappa_{j_0}(T,B,s,r)&=\int_{\Tube(B,s,r)}\tau_*T\wedge \pi^*(\omega^{j_0})\wedge\alpha^{q_0}
  = I_{q_0,0,j_0,0}(T,s,r)\\&+\sum_{j_1,j'_1} {q_0 \choose j_1}
  {q_0-j_1 \choose j'_1}(-1)^{q_0-j_1-j'_1}(-c_2)^{j'_1}
 I_{j_1, 0, q_0+j_0-j_1-j'_1,q_0-j_1-j'_1}(T,s,r).
\end{split}
\end{equation}
Recalling that $q_0=k-p-j_0,$  we set
\begin{equation}\label{e:Cal-I}
\begin{split}
 \Ic_{j_0}(T,s,r)&:=\kappa_{j_0}(T, B,s, r,\tau,h)- \sum_{\lowm\leq j_1<q_0} {q_0\choose j_1} (-1)^{q_0-j_1}c_1^{q_0-j_1} \Kc_{j_1,k-p-j_1}(T, B,s, r)\\
 &-\Kc_{j_0,k-p-j_0}(T, B,s,r).
 \end{split}
\end{equation}
 We rewrite   \eqref{e:P-Lc-bullet_finite(5)} as 
 \begin{equation}\label{e:P-Lc-bullet_finite(6)}
 \Ic_{j_0}(T,s,r)
 =\Ic_1+\Ic_2+\Ic_3,
\end{equation}
where
\begin{equation*}
 \begin{split}
 \Ic_1&:=\sum_{j'_1,j_1}{q_0 \choose j_1}
  {q_0-j_1 \choose j'_1}(-1)^{q_0-j_1-j'_1}(-c_2)^{j'_1}  
 \cdot I^\hash_{j_1, 0, q_0+j_0-j_1-j'_1,q_0-j_1-j'_1}(T,r), \\
 \Ic_2&:=    \sum_{\lowm\leq j_1\leq q_0} {q_0\choose j_1} (-1)^{q_0-j_1}c_1^{q_0-j_1}   \big(   I_{k-p-j_1,0,j_1,0}(T,s,r)-I^\hash_{k-p-j_1,0,j_1,0}(T,s,r)\big), \\
 \Ic_3&:=  \sum_{j'_1,j_1} {q_0 \choose j_1}
  {q_0-j_1 \choose j'_1}(-1)^{q_0-j_1-j'_1}(-c_2)^{j'_1}\\
  &\cdot\big(I_{j_1, 0, q_0+j_0-j_1-j'_1,q_0-j_1-j'_1}(T,r)-I^\hash_{j_1, 0, q_0+j_0-j_1-j'_1,q_0-j_1-j'_1}(T,r)\big)  .
\end{split}
\end{equation*}
Here the  sums $\Ic_1$ and $\Ic_3$  are  taken over all $(j_1,j'_1)$   such that   $0\leq j_1,$ $0<j'_1$ and
$j_1+j'_1\leq  q_0.$ 

Consider an arbitrary term
$I^\hash_{j_1, 0, q_0+j_0-j_1-j'_1,q_0-j_1-j'_1}(T,r)$ in the sum on the  expression of $\Ic_1.$
Since $j_1+(q_0+j_0-j_1-j'_1)=q_0+j_0-j'_1=k-p-j'_1<k-p$ as $j'_1>0,$  the term is $\lesssim \Kc^-(T,s,r). $  
Consequently, we get that
\begin{equation}\label{e:P-Lc-bullet_finite(6bis)}
  |\Ic_1|\leq  c  \Kc^-(T,s,r). 
\end{equation}
Applying Lemma \ref{L:spec-wedge}  to  each  difference term  in $\Ic_2$ and $\Ic_3$ 
yields that
\begin{equation} \label{e:P-Lc-bullet_finite(7)}
|I_\bfi(r)- I^\hash_\bfi(s,r)|^2 \leq c\big(\sum_{\bfi'} I^\hash_{\bfi'}(s-c_0s^2,r+c_0r^2)\big)\big ( \sum_{\bfi''} I^\hash_{\bfi''}(s-c_0s^2,r+c_0r^2) \big).  
\end{equation}
Here, on the LHS  $\bfi=(i_1,i_2,i_3,i_4)$  is  either $(q_0,0,j_0,0)$ or $(j_1,0, q_0+j_0-j_1-j'_1,q_0-j_1-j'_1)$  with $j_1,j'_1$ as above, 
and on the RHS:
\begin{itemize} \item[$\bullet$] the first sum  is taken over a finite number of multi-indices    $\bfi'=(i'_1,i'_2,i'_3,i'_4)$ as above  such that  $i'_1\leq  i_1$  and $i'_2\geq i_2;$ 
\item  the second sum   is taken over  a finite number of multi-indices $\bfi''=(i''_1,i''_2,i''_3,i''_4)$ as above   such that   either  ($i''_1< i_1$)
or ($i''_1=i_1$ and $i''_2\geq {1\over 4}+i_2$) or ($i''_1=i_1$ and $i''_3<i_3$).
\end{itemize}
 Observe that   $c_1-c_2\varphi\approx  1$ on $\Tube(B,\bfr)$  and that $\varphi \lesssim r^2\lesssim 1$ on $\Tube(B,r)$ for $0<r\leq\bfr.$
 Therefore, $I_{i_1,i_2,i_3,i_4}(T,s,r)\leq c I_{i_1,0,i_3,0}(T,s,r)$ for a constant $c>0$ independent of $T$ and $0<r\leq \bfr.$
 Consequently, the first  sum on the RHS of \eqref{e:P-Lc-bullet_finite(7)} is  bounded from above by a constant times  $\Kc_{q_0}(T,s-c_0s^2,r+c_0r^2),$
 whereas  the second sum  is bounded from above by a constant times $\Kc^-(T,s-c_0s^2,r+c_0r^2)+ r^{1\over 2} \Kc_{q_0}(T,s-c_0s^2,r+c_0r^2).$
 In fact the factor  $r^{1\over 2}$ comes  from  $\varphi^{1\over 4}$ because $\varphi\lesssim r^2$ on $\Tube(B,r).$    
Consequently,   we infer from \eqref{e:Cal-I}--\eqref{e:P-Lc-bullet_finite(6)}--\eqref{e:P-Lc-bullet_finite(6bis)} and \eqref{e:P-Lc-bullet_finite(7)} that there is a constant $c>0$ such that
\begin{equation*}
 \begin{split}
|\Ic_{j_0}(T,s,r)|
 &\leq    c\Kc^-(T,s,r) \\
 &+c\sqrt{\Kc_{q_0}(T,s-c_0s^2,r+c_0r^2)} \sqrt{ \Kc^-(T,s-c_0s^2,r+c_0r^2)+r^{1\over 2} \Kc_{q_0}(T,s-c_0s^2,r+c_0r^2) }.
\end{split}
 \end{equation*}
By Proposition \ref{P:main},  $\Kc^-(T,s,r)\lesssim r^2$ and $\Kc_q(T,s,r)\lesssim 1.$ Hence, 
\begin{equation}\label{e:P-Lc-bullet_finite(8)}
 |\Ic_{j_0}(T,s,r)|\leq    cr^{1\over 4}
 .
\end{equation}
We  will prove by decreasing induction on $j_0:\  \lowm\leq j_0\leq \upm$ that there is a constant $c>0$ such that 
\begin{equation}\label{e:induction}   c^{-1}\big(\sum_{j=\lowm}^{j_0} \Kc_{j,k-p-j}(T, B,s,r) +r^{1\over 4} \big)\leq  \sum_{j=\lowm}^{j_0} |  \kappa_j(T, B,s, r,\tau,h)|+ r^{1\over 4}  \leq c\big(\sum_{j=\lowm}^{j_0} \Kc_{j,k-p-j}(T, B,s,r) +r^{1\over 4} \big).
\end{equation}
Applying \eqref{e:Cal-I} and \eqref{e:P-Lc-bullet_finite(8)} for $j=\upm$ yields that
\begin{equation*}
|\kappa_{\lowm}(T, B,s, r,\tau,h)-  \Kc_{\lowm,k-p-\lowm}(T, B,s, r)|\leq    cr^{1\over 4}.
\end{equation*}
This  proves \eqref{e:induction} for $j_0=\upm.$

Suppose \eqref{e:induction} for $j_0+1.$ We need to show  \eqref{e:induction} for $j_0.$
 Putting  \eqref{e:Cal-I} and \eqref{e:P-Lc-bullet_finite(8)} together, we get that
 \begin{eqnarray*}
 |\kappa_{j_0}(T, B,s, r,\tau,h)|+r^{1\over 4}&\leq& c( \sum_{\lowm\leq j_1<q_0} {q_0\choose j_1} c_1^{q_0-j_1} \Kc_{j_1,k-p-j_1}(T, B,s, r)+\Kc_{j_0,k-p-j_0}(T, B,s,r)+r^{1\over 4}),\\
 \Kc_{j_0,k-p-j_0}(T, B,s,r)+r^{1\over 4}&\leq& c( \sum_{\lowm\leq j_1<q_0} {q_0\choose j_1} c_1^{q_0-j_1} \Kc_{j_1,k-p-j_1}(T, B,s, r)+|\kappa_{j_0}(T, B,s, r,\tau,h)|+r^{1\over 4}).
 \end{eqnarray*}
 This, combined with  the  induction hypothesis (that is, \eqref{e:induction} for $j_0+1$),
 completes the proof of \eqref{e:induction} for $j_0.$  The proof of assertion (1) is thereby completed.

The proof of assertion (2) is  similar.
 \endproof
 
\begin{proposition}\label{P:main-1}
 We keep the   hypothesis and the notation  of  Theorem \ref{T:main_1}.
    \begin{itemize}
  \item [(a)] If condition  (a)  of Theorem \ref{T:main_1} is satisfied,  then  $$\int_0^{r_0} {\sqrt{\Kc_{j,k-p-j}(T, B,r/2,r)}\over  r}dr <\infty\qquad \text{for all}\qquad \lowm\leq j\leq \upm;$$
  \item [(b)] If condition  (b)  of Theorem \ref{T:main_1} is satisfied,  then $$ \int_0^{r_0} {\Kc_{j,k-p-j}(T,B,0,r)\over  r}dr <\infty\qquad \text{for all}\qquad \lowm\leq j\leq \upm.$$
  \end{itemize}
 \end{proposition}
 \proof
 First we will prove assertion (a). So  we  assume that condition  (a)  of Theorem \ref{T:main_1} is satisfied.
 Fix an index $j$ with  $\lowm\leq j\leq \upm.$  
 Recall from  the  hypothesis that  $T=T^+-T^-$ with $T^\pm\in \CL^{p;1,1}(B).$
 Let $(T^\pm_n)_{n=0}^\infty$ be an approximating    sequence of $T^\pm$ according  to  Definition \ref{D:Class}.
 Set $T_n:=T^+_n-T^-_n$ for $n\in\N.$
 By Lelong-Jensen formula for closed currents (see Theorem \ref{T:Lelong-Jensen-closed}), 
 for  all $r\in  (0,\bfr]$ except for  a  countable set of   values, we have that  
\begin{equation}\label{e:Lelong-Jensen_s=r/2}
 \nu_j(T,B,r,\tau,h)-  \nu_j(T,B,r/2,\tau,h)
 =\lim\limits_{n\to\infty} \Vc(T_n,r)+   \kappa_j(T,B,r/2,r,\tau,h).
\end{equation} 
Here the vertical boundary term  $\Vc(S,r)$ for  a  continuous form $S$  is  given by
\begin{multline*}
 \Vc(S,r):=
 {1\over r^{2(k-p-j)}} \int_{\partial_\ver \Tube(B,r)}\dc\varphi\wedge S\wedge \beta^{k-p-j-1}\\
 -{1\over (r/2)^{2(k-p-j)}} \int_{\partial_\ver \Tube(B,r/2)}\dc\varphi\wedge S\wedge \beta^{k-p-j-1}-\int_{\partial_\ver \Tube(B,r/2,r)}\dc\log\varphi\wedge S\wedge \alpha^{k-p-j-1}.
\end{multline*}
Since  $T_n=T^+_n-T^-_n$ and $(T^\pm_n)$ is an approximating sequence of $T^\pm\in\CL^{p;1,1}(B)$ in the sense of Definition \ref{D:Class}, there is a  constant $c$ depending only on $T$  such  that 
\begin{equation*}  \big| \Vc(T_n,r)\big|\leq c r
 \qquad\text{for}\qquad 0<r\leq\bfr,\qquad n\geq 1.
\end{equation*}
On the other hand, by  
Proposition \ref{P:main}  applied  to $s=r/2,$ we obtain that
\begin{equation*}     \sum_{j=\lowm}^\upm\Kc_{j,k-p-j}(T, B,r/2,r)\leq c\sum_{j=\lowm}^\upm |  \kappa_j(T, B,r/2,r,\tau,h)|+ cr^{1/4}\qquad\text{for}\qquad 0<r\leq \bfr.
\end{equation*}
Inserting the last two inequalities in  equality \eqref{e:Lelong-Jensen_s=r/2}, we infer that
\begin{equation*} 
 \sum_{j=\lowm}^\upm\Kc_{j,k-p-j}(T, B,r/2,r)\leq c  \sum_{j=\lowm}^\upm \big| \nu_j(T,B,r,\tau,h)-  \nu_j(T,B,r/2,\tau,h)\big|  +cr^{1\over 4}\qquad\text{for}\qquad 0<r\leq \bfr.
\end{equation*} 
This, combined with condition (a) and  the  elementary equality $\sqrt{x+y}\leq \sqrt{x}+\sqrt{y}$ for $x,y\geq 0,$ implies that
  \begin{eqnarray*} &&\int_0^{r_0} {\sqrt{\sum_{j=\lowm}^\upm\Kc_{j,k-p-j}(T, B,r/2,r)}\over  r}dr\\
&\leq  & \int_0^{r_0} {\sqrt{ \sum_{j=\lowm}^\upm |\nu_j(T,B,r,\tau,h)-  \nu_j(T,B,r/2,\tau,h)|+cr^{1\over 4        }}\over  r}dr
\\
  &\leq&\sum_{j=\lowm}^\upm\int_0^{r_0} {\sqrt{  |\nu_j(T,B,r,\tau,h)-  \nu_j(T,B,r/2,\tau,h)|       }\over  r}dr + c\int_0^{r_0}{r^{1\over 4}         \over  r}dr.
  \end{eqnarray*}
Since the first  term  in the last  line  is    finite by condition (a) and the second term  in the last line is clearly finite, assertion (a) follows.  

 It remains to prove  assertion (b). So  we  assume that condition  (b)  of Theorem \ref{T:main_1} is satisfied.
 Fix an index $j$ with  $\lowm\leq j\leq \upm.$  By Lelong-Jensen formula for closed currents (see Theorem \ref{T:Lelong-Jensen-closed}), 
 for  all $0<s<r\in  (0,\bfr]$ except for  a  countable set of   values, we have that  
\begin{equation}\label{e:Lelong-Jensen_s,r}
 \nu_j(T,B,r,\tau,h)-  \nu_j(T,B,s,\tau,h)
 =\lim\limits_{n\to\infty} \Vc(T_n,s,r)+   \kappa_j(T,B,s,r,\tau,h).
\end{equation} 
Since  $T$ belongs to the class $\CL^{p;1,1}(B),$ there is a  constant $c>0$ depending only on $T$  such  that 
\begin{equation*}  \big| \Vc(T_n,s,r)\big|\leq c r
 \qquad\text{for}\qquad 0<r\leq\bfr,\qquad n\geq 1.
\end{equation*}
By Theorem \ref{T:Nguyen} (1), $\lim_{s\to 0+} \nu_j(T,B,s,\tau,h)=\nu_j(T,B,h).$ 
Morever,   $\lim_{s\to 0+} \kappa_j(T,B,s,t,\tau,h)=\kappa^\bullet_j(T,B,r,h).$  Therefore,
 letting $s$ tend to $0,$  we  infer that from the last estimate and \eqref{e:Lelong-Jensen_s,r} that
\begin{equation*} 
\big| (\nu_j(T,B,r,\tau,h)-  \nu_j(T,B,h)) - \kappa^\bullet_j(T,B,r,\tau,h)\big|\leq  c r
 \qquad\text{for}\qquad 0<r\leq\bfr .
\end{equation*} 
On the other hand, by  
Proposition \ref{P:main-bis}  applied  to $s=0,$ we obtain that
\begin{equation*}   \sum_{j=\lowm}^\upm\Kc_{j,k-p-j}(T, B,0,r)   \leq c\sum_{j=\lowm}^\upm \big |\kappa^\bullet_j(T, B,r,\tau,h) \big|+  cr^{1/4}\qquad\text{for}\qquad 0<r\leq \bfr.
\end{equation*}
We infer from  the last two inequalities that
\begin{equation*} 
  \sum_{j=\lowm}^\upm\Kc_{j,k-p-j}(T, B,0,r)   \leq  \sum_{j=\lowm}^\upm  \big| \nu_j(T,B,r,\tau,h)-  \nu_j(T,B,h) \big|+ cr^{1\over 4}\qquad\text{for}\qquad 0<r\leq \bfr.
\end{equation*} 
This, combined with condition (b), implies that
  \begin{equation*} \sum_{j=\lowm}^\upm \int_0^{r_0} {\Kc_{j,k-p-j}(T, B,0,r)\over  r}dr\leq   \sum_{j=\lowm}^\upm \int_0^{r_0} {  |\nu_j(T,B,r,\tau,h)-  \nu_j(T,B,h)|       \over  r}dr +c\int_0^{r_0} r^{ -3\over4}dr  .
  \end{equation*}
Since the first  term  in the last  line  is    finite by condition (b) and the second term  in the last line is clearly finite, assertion (b) follows.  
 \endproof
  
\begin{proposition}\label{P:main-2}
 We keep the   hypothesis and the notation  of  Theorem \ref{T:main_2} and   let $i\in I.$
    \begin{itemize}
  \item [(a)] If condition  (a-i)  of Theorem \ref{T:main_2} is satisfied,  then  $$\int_0^{r_0} {\sqrt{\Kc_{j,k-p-j}(T, B_i\cap B,r/2,r)}\over  r}dr <\infty\qquad \text{for all}\qquad \lowm\leq j\leq \upm;$$
  \item [(b)] If condition  (b-i)  of Theorem \ref{T:main_2} is satisfied,  then $$ \int_0^{r_0} {\Kc_{j,k-p-j}(T,B_i\cap B,0,r)\over  r}dr <\infty\qquad \text{for all}\qquad \lowm\leq j\leq \upm.$$
  \end{itemize}
 \end{proposition}
 \proof
 First we prove  assertion (a).
 Condition (a-i), combined  with the  elementary equality $\sqrt{x+y}\leq \sqrt{x}+\sqrt{y}$ for $x,y\geq 0,$ implies that
  \begin{equation*} \int_0^{r_0} {\sqrt{\sum_{j=\lowm}^\upm |  \kappa_j(T,B_i\cap  B,r/2, r,\tau,h)|}\over r}dr 
   \leq \sum_{j=\lowm}^\upm\int_0^{r_0} {\sqrt{|  \kappa_j(T,B_i\cap  B,r/2, r,\tau,h)|}\over r}dr 
   <\infty.
  \end{equation*}
By Proposition \ref{P:main-bis}  there is a constant $c>0$  such that for every
   $T\in \widetilde\CL^{p;1,1} (\bfU,\bfW)$ and $0<r\leq \bfr,$  
\begin{equation*}\sum_{j=\lowm}^\upm \Kc_{j,k-p-j}(T, B_i\cap B,r/2,r) \leq c\big( \sum_{j=\lowm}^\upm |  \kappa_j(T,B_i\cap  B,r/2, r,\tau,h)|+ r^{1\over 4} \big).
\end{equation*}
This, coupled with the previous  estimate, implies that
\begin{equation*} \int_0^{r_0} {
   \sqrt{    \sum_{j=\lowm}^\upm    \Kc_{j,k-p-j}(T, B_i\cap B,r/2,r)  }      \over r}dr 
   \lesssim  \sum_{j=\lowm}^\upm\int_0^{r_0} {\sqrt{|  \kappa_j(T,B_i\cap  B,r/2, r,\tau,h)|}\over r}dr
   +\int_0^{r_0}{r^{-7\over 8}}dr
   .
  \end{equation*}
  This  proves assertion (a).
  
  We turn to the proof of assertion (b). By Proposition \ref{P:main-bis}  there is a constant $c>0$  such that for every
   $T\in \widetilde\CL^{p;1,1} (\bfU,\bfW)$ and $0<r\leq \bfr,$  
\begin{equation*}\sum_{j=\lowm}^\upm \Kc_{j,k-p-j}(T, B_i\cap B,0,r) \leq c\big( \sum_{j=\lowm}^\upm |  \kappa^\bullet_j(T,B_i\cap  B, r,\tau,h)|+ r^{1\over 4} \big).
\end{equation*}
This,  combined with  condition (b-i), implies that
\begin{equation*} \int_0^{r_0} {
      \sum_{j=\lowm}^\upm    \Kc_{j,k-p-j}(T, B_i\cap B,0,r)        \over r}dr 
   \lesssim  \sum_{j=\lowm}^\upm\int_0^{r_0} {|  \kappa_j(T,B_i\cap  B,r/2, r,\tau,h)|\over r}dr
   +\int_0^{r_0}{r^{-3\over 4}}dr
   .
  \end{equation*}
  Hence, assertion (b) follows.

\endproof  
\section{Local  core estimates} \label{S:Local-estimates}
 
 This  section  infers  from the    global estimates  established in  Section \ref{S:Global-estimates} corresponding local     growth behavior of appropriate local mass indicators.

\subsection{Analysis  in local coordinates}
Since $V_0\Subset  V,$ we only need to  prove  a local result near a given point $y_0\in V_0.$
We use the coordinates $(z,w)\in\C^{k-l}\times \C^l$  around a neighborhood $U$ of $y_0$  such that $y_0=0$  in  these coordinates. 
We may assume that  $U$ has the form $U=U'\times U'',$ where $U'$ (resp. $U'')$ is an  open neighborhood of $0'$ in $\C^{k-l}$ (resp. of  $0''$ in $\C^l$), 
and  $V=\{z=0\}\simeq U''.$ Moreover,  we  may assume  that $U''=(2\D)^l.$  
Consider  the trivial  vector bundle $\pi:\ \E \to  U''$ with  $\E\simeq  \C^{k-l}\times U''.$ 

Let $T$ be a $(p,p)$-current of order $0$  defined on  an open neighborhood $U$ of $0$ in $\C^k.$
We use the coordinates $y=(z,w)\in\C^{k-l}\times \C^l.$ 
We may assume that  $U$ has the form $U=U'\times U'',$ where $U'$ (resp. $U'')$ are open neighborhood of $0'$ in $\C^{k-l}$ of  ($0''$ in $\C^l$).
Let $V=\{z=0\}=U''$ and let $B=B_w\Subset U''$ be a domain with  piecewise  $\Cc^2$-smooth boundary    and $\bfr>0$  such that  $\{\|z\|<\bfr\}\times B\Subset U.$ 
Consider  the trivial  vector bundle $\pi:\ \E \to  U''$ with  $\E\simeq  \C^{k-l}\times U''$ endowed with the canonical Euclidean  metric $h$ on fibers. For $\lambda\in\C^*,$  let $a_\lambda:\ \E\to \E$ be the multiplication by  $\lambda$
on fibers, that is, 
$a_\lambda(z,w):=(\lambda z,w)$ for $(z,w)\in \E.$
The admissible map $\tau$ in this setting is simply the identity $\id.$

Consider the positive closed $(1,1)$-forms
\begin{equation}\label{e:alpha-beta-Upsilon-local}\beta=\omega_z:=\ddc \|z\|^2\quad\text{and}\quad  \omega=\omega_w:=\ddc\|w\|^2\quad\text{and}\quad \alpha:=\ddc \log\|z\|^2.
\end{equation}
In fact,  by \eqref{e:alpha-beta-spec} and  \eqref{e:varphi-new-exp}, we have $\alpha:=\ddc \log\|A(w)z\|^2.$ By the above choice of $h,$
$A(w)$ is the identity matrix in $\C^{k-l}.$
Let $\lowm\leq j\leq \upm.$ For  $0<r<\bfr,$ consider the quantity 
\begin{equation}\label{e:Lelong-numbers-local}
 \nu_j(T,B,r):=  {1\over r^{2(k-p-j)}}\int_{\|z\|<r,\ w\in B} T\wedge \omega_w^j \wedge \omega_z^{k-p-j}.
\end{equation}
 For  $0<s<r\leq \bfr,$   consider 
\begin{equation}\label{e:Lelong-corona-numbers-local}
 \kappa_j(T,B,s,r):=   \int_{s<\|z\|<r,\ w\in B}  T\wedge \omega_w^j \wedge \alpha^{k-p-j}. 
\end{equation} 
Let $0<r\leq \bfr.$   Consider
\begin{equation}\label{e:Lelong-log-bullet-numbers-local}
 \kappa^\bullet_j(T,B,r):= \limsup\limits_{s\to0+}   \kappa_j(T,B,s,r). 
\end{equation} 
We  also consider
\begin{equation}\label{e:Lelong-log-numbers-local}
 \kappa_j(T,B,r):=  \int_{\|z\|<r,\ w\in B}  T\wedge \omega_w^j \wedge \alpha^{k-p-j},  
\end{equation}
provided that the right hand  side makes sense.

  Following the analysis in local coordinates  in \cite[Section 5.2]{Nguyen21},
let $\pi_\FS:\ \C^{k-l}\setminus \{0\}\to  \P^{k-l-1},$  $z\mapsto  \pi_\FS(z):=[z]$   be the canonical projection, and 
let $\omega_\FS$ be the  Fubini-Study  form  on  $\P^{k-l-1}.$  So
\begin{equation}\label{e:FS}
\pi^*_\FS  (\omega_\FS([z]))=\ddc (\log{\|z\|^2})\qquad\text{for}\qquad  z\in\C^{k-l}\setminus \{0\}.
\end{equation}
We place ourselves on  an open set $U$ of $\C^{k-l}$ defined by $z_{k-l}\not=0.$
We   may assume without loss of generality that
\begin{equation}\label{e:max-coordinate} 2|z_{k-l}| > \max\limits_{1\leq j\leq k-l}|z_j|.
\end{equation}
and use the ``projective'' coordinates
\begin{equation}\label{e:homogeneous-coordinates}
\zeta_1:={z_1\over z_{k-l}},\ldots, \zeta_{k-l-1}:={z_{k-l-1}\over z_{k-l}},\quad \zeta_{k-l}=z_{k-l}.
\end{equation}
Therefore,  the coordinates  $\zeta=(\zeta_1,\ldots,\zeta_{k-l})=(\zeta',\zeta_{k-l})$ vary in  an open set $ \widetilde U\subset \C^{k-l}$ which  depends only on $U.$  With these new coordinates, the form $\omega_\FS([z])$  can be  rewritten as  
\begin{equation}\label{e:FS-zeta} \omega_\FS([z])= \ddc \log{ (1+|\zeta_1|^2+\cdots+|\zeta_{k-l-1}|^2)},
\end{equation}
and a direct computation shows that 
\begin{equation}\label{e:beta'-zeta'}
\omega_\FS([z])\approx  (1+\|\zeta'\|^2)^{-2}\beta'(\zeta'),\quad\text{where}\quad \beta'(\zeta'):=\ddc (|\zeta_1|^2+\cdots+|\zeta_{k-l-1}|^2).
\end{equation}
Since  $|\zeta_j|<2$ for $1\leq j\leq k-l-1$ by  \eqref{e:max-coordinate}, it follows from  \eqref{e:beta'-zeta'}   that
\begin{equation}\label{e:omega_FS-vs-beta'}
\omega_\FS([z])\approx \beta'(\zeta') .
\end{equation}
  
  Let $S$ be  a differential  form (resp.  a  current) defined on $\Tube(B,r)\subset \E$ 
  for some $0<r\leq\bfr.$
   So we can write  in a  local representation of $S$ in coordinates $y=(z,w)\in\C^{k-l}\times\C^l: $
   \begin{equation}\label{e:S_IJKL} S=\sum_{M=(I,J;K,L)} S_Mdz_I\wedge d\bar z_J\wedge dw_K\wedge d\bar w_L,
   \end{equation}
where the $S_M=S_{I,J;K,L}(z,w)$  are  the component  functions  (resp.  component  distributions), and the sum is taken over   $M=( I,J;K,L)$ with $I,J\subset\{1,\ldots,k-l\}$ and $K,L\subset \{1,\ldots,l\}.$

For  $M=( I,J;K,L)$ as above, we also write $dy_M$ instead of $dz_I\wedge d\bar z_J\wedge dw_K\wedge d\bar w_L.$
 

 \subsection{Main  estimates}

 Let $\lambda\in\C^*$ and set $t:=\lambda^{-1}.$ So in the coordinates $(\zeta',\zeta_{k-l},w),$ the 
 dilation  $A_t$   has  the form
 $$
 A_t:\qquad  (\zeta',\zeta_{k-l},w)\mapsto  (\zeta',t \zeta_{k-l},w).
 $$
 So 
 the coefficients of $(A_\lambda)_*T= (A_t)^*T$ are given  by
 \begin{equation}\label{e:coeff-A-lambda_T}
  T_{I;J;K,L}^t (\zeta',\zeta_{k-l},w)=\begin{cases} T_{I,J;K,L} (\zeta', t\zeta_{k-l},w), & \text{if}\  k-l\not\in I,\ k-l\not\in J;\\
  t T_{I,J;K,L} (\zeta', t\zeta_{k-l},w), & \text{if}\  k-l\in I,\ k-l\not\in J;\\
       \bar t T_{I,J;K,L} (\zeta', t\zeta_{k-l},w), & \text{if}\  k-l\not\in I,\ k-l\in J;\\  
        |t|^2 T_{I,J;K,L} (\zeta', t\zeta_{k-l},w), & \text{if}\  k-l\in I,\ k-l \in J.
                       \end{cases}
 \end{equation}
Consider the functions $\kappa,\ \kappa^\bullet:\ (0,\bfr]\to \R^+$ given by
\begin{equation}\label{e:two-kappa}
 \kappa(r):=\sum_{j=\lowm}^\upm\kappa_j(T,B,r/2,r)\quad \text{and}\quad
  \kappa^\bullet(r):=\sum_{j=\lowm}^\upm\kappa^\bullet_j(T,B,r)\quad \text{for}\quad  r\in (0,\bfr].
\end{equation}
Recall  from  \cite[Theorems  3.5, 3.6 and 3.7] {Nguyen21} that  $\kappa(r) ,\kappa^\bullet(r)\to 0$
as $r\to 0.$
If condition  (a)  of Theorem \ref{T:main_1} is satisfied,  then by Proposition \ref{P:main-1},
  $$\int_0^{r_0} {\sqrt{\Kc_{j,k-p-j}(T, B,r/2,r)}\over  r}dr <\infty\qquad \text{for all}\qquad \lowm\leq j\leq \upm;$$
and hence  we infer that
$$
\int_0^{r_0}{\sqrt{\kappa(r)}\over r}dr <\infty.
$$
If condition  (b)  of Theorem \ref{T:main_1} is satisfied,  then by Proposition \ref{P:main-1},
  $$\int_0^{r_0} {\Kc_{j,k-p-j}(T, B,0,r)\over  r}dr <\infty\qquad \text{for all}\qquad \lowm\leq j\leq \upm;$$
and hence  we infer that
$$
\int_0^{r_0}{\kappa^\bullet(r)\over r}dr <\infty.
$$
On the other hand, we deduce from  Theorem \ref{T:basic-estimate} that  
\begin{equation}\label{e:local-initial-est-1}
 \int_{\Tube(B,r/2,r)} (A_\lambda)_*T\wedge (\pi^*\omega)^j\wedge \alpha^{k-p-j}=\kappa_j(|t|r)\leq c.
\end{equation}
We also have
\begin{equation}\label{e:local-initial-est-2}
 \int_{\Tube(B,r)} (A_\lambda)_*T\wedge (\pi^*\omega)^j\wedge \beta^{k-p-j}\leq cr^{2(k-p-j)}.
\end{equation}
The next lemma allows to bound the coefficients in \eqref{e:coeff-A-lambda_T} by the functions in \eqref{e:two-kappa}.
 \begin{lemma}\label{L:coeff-Tt}
  There are constants  $c_1,c_2,c_3>0 $ such that  the measures   $T_{I;J;K,L}^t$ satisfy the following    inequalities
  \begin{equation*}
  \int_{\widetilde U}|T_{I;J;K,L}^t|\leq \begin{cases} c_1, & \text{for all }\  I, J,K,L;\\
  c_2\kappa (|t|) , & \text{if}\  k-l\in I,\ k-l\in J;\\
       c_3\sqrt{\kappa (|t|)} , & \text{if}\  k-l\in I\ \text{or}\  k-l\in J  .
                       \end{cases}
 \end{equation*}
 \end{lemma}
\proof
Let $y^{(0)}=(z^{(0)},w^{(0)})\in(\C^{k-l}\setminus\{0\})\times \C_l.$
We may assume  that  
\begin{equation}\label{e:max-coordinate-bis} 2|z^{(0)}_{k-l}| > \max\limits_{1\leq j\leq k-l}|z^{(0)}_j|.
\end{equation}
Without loss of generality we may assume that ${1\over 2}\leq |z^{(0)}_{k-l}|\leq 1.$
So  we  infer that
\begin{equation*}
 \alpha\geq  c\beta' \quad\text{on}\qquad\widetilde U\setminus V.
\end{equation*}
Using  estimate \eqref{e:local-initial-est-2}, we  deduce  that $\int_U|T_{I;J;K,L}^t|\leq c_1$
for all $I,J,K,L.$ 

On the other hand,
Theorem \ref{T:Nguyen} implies that
 \begin{equation*}
 \int_{\Tube(B,r/2,r)} (A_t)^*T\wedge (\pi^*\omega+\beta')^{k-p}\leq \kappa(|t|r)\leq c.
\end{equation*}
Hence, we infer that
\begin{equation*} 
 \int_{\Tube(B,1/2,1)} (A_t)^*T\wedge (\pi^*\omega+\beta')^{k-p}\leq \kappa(|t|)\leq c.
\end{equation*}
Consequently, we get that
\begin{equation} \label{e:local-initial-est-bis}
 \int_{\Tube(B,1/2,1)} \sum_{I\ni k-l} T^t_{I,I;K,K}\leq \kappa(|t|)\leq c.
\end{equation}
Applying  Proposition  \ref{P:Demailly} to estimate \eqref{e:local-initial-est-bis} yields that
 $ \int_{\widetilde U}|T_{I;J;K,L}^t|< c_2\kappa (|t|) $ if $  k-l\in I $ and $k-l\in J.$
 
 It remains  to treat the case  where $ k-l\in I$ and $  k-l\not\in J.$
 Choose $\lambda_1=\ldots=\lambda_{k-l-1}=1$ and $\lambda_{k-l+1}=\ldots=\lambda_k=1.$
 
Applying  Proposition  \ref{P:Demailly} to estimate \eqref{e:local-initial-est-bis} yields that
\begin{equation*}
\lambda_{k-l}  \int_{\widetilde U}|T_{I;J;K,L}^t|\leq  \int_{\widetilde U\setminus V}  \sum_{M=(I,I;K,K):\ k-l\not\in I}|T_{M}^t|
 +\lambda_{k-l}^2\int_{\widetilde U}  \sum_{M=(I,I;K,K):\ k-l\in I}|T_{M}^t|\le (c_0+c_1|\lambda_{k-l}|^2\kappa(|t|). 
\end{equation*}
The desired  last estimate follows by choosing $\lambda_{k-l}:=\kappa(|t|)^{-1\over 2}.$
\endproof

\section{Proof of the criteria}\label{S:Proofs}

Having  good local estimates obtained  in Section \ref{S:Local-estimates} in hands,  we employ in  this  section the limiting argument of Blel-Demailly-Mouzali \cite{BlelDemaillyMouzali} in order to get the 
desired convergence of $(A_\lambda)_*(\tau_*T)$ as $\lambda$ tends to infinity.

More  specifically,  we prove that  $T^t_{I,J;K,L}$ converges weakly on $U$ as $t$ tends to $0.$
By Lemma \ref{L:coeff-Tt}, the coefficient $T^t_{I,J;K,L}$ tends to $0$ in mass when either $I$ or  $J$ contains $k-l.$
It remains  to show that  $T^t_{I,J;K,L}$ converges weakly   when neither $I$ nor  $J$ contains $k-l.$
 
 Let $\phi$ be  a smooth function compactly supported in $\E|_B.$  Let  $I,J$  such that $k-l\not\in I,$ $k-l\not\in J.$
 Consider  the  function
 \begin{equation}\label{e:f_I,J;K,L}
  f_{I,J;K,L}(t):=\int_{\widetilde U} T^t_{I,J;K,L}(\zeta,w)\phi(\zeta,w)d\Leb(\zeta,w)=\int_U T_{I,J;K,L}(\zeta',t\zeta_{k-l},w)\phi(\zeta,w)d\Leb(\zeta,w).
 \end{equation}
The  functions $f_{I,J;K,L}$'s   are of $\Cc^\infty$ class on the  punctured  disc  $\{t\in\C:  0<|t|<1\}$ and is  bounded in a  neighborhood of $0.$
The  problem  is to check if    $f_{I,J;K,L}(t)$  admits a limit    when  $|t|$ tends to $0.$   To this  end  we  estimate  the  derivatives 
${\partial f_{I,J;K,L}\over \partial t},$ ${\partial f_{I,J;K,L}\over \partial \bar t},  $ ${\partial^2 f_{I,J;K,L}\over \partial t\partial \bar t}$ in a neighborhood of $0.$ 
\begin{lemma}\label{L:f_I,J;K,L}
 There are constants  $c_1,c_2>0$ such that for   $0<|t|\leq \bfr,$  and for 
 $I,J\subset\{1,\ldots,k-l\}$ and $K,L\subset \{1,\ldots,l\}$   satisfying $k-l\not\in I$ and $k-l\not\in J,$  it holds that
 \begin{equation*}
  \big| {\partial  f_{I,J;K,L} \over \dbar t}\big|\leq c_1{\sqrt{\kappa(|t|)}\over |t|}\quad\text{and}\big| {\partial  f_{I,J;K,L} \over \partial t}\big|\leq c_1{\sqrt{\kappa(|t|)}\over |t|}\quad\text{and}\quad  \big| {\partial^2  f_{I,J;K,L} \over
 \partial t \partial\bar  t}\big|\leq c_2 {\kappa(|t|)\over |t|^2}
 \end{equation*}
\end{lemma}
\proof
By a differentiation  under the  integral sign we get that
 \begin{equation}\label{e:deri_f_I,J;K,L_prem}
 {\partial  f_{I,J;K,L} \over \partial t}=\int_U \zeta_{k-l} {\partial T_{I,J;K,L}(\zeta',t\zeta_{k-l},w)\over  \partial \zeta_{k-l}}\phi(\zeta,w)d\Leb(\zeta,w).
 \end{equation}
 In what follows  we  adopt the  following notation
 \begin{eqnarray*}
  I(q)&:=& (I\setminus \{i_{q}\})\cup\{k-l\}\qquad \text{for}\quad I=(i_1,\ldots,i_{|I|});\\
 J(q)&:=& (J\setminus \{j_{q}\})\cup\{k-l\}\qquad \text{for}\quad J=(j_1,\ldots,j_{|J|});\\
 K(q)&:=& K\setminus \{k_{q}\}\qquad\text{for}\qquad  K=(k_1,\ldots,k_{|K|});\\
  L(q)&:=& L\setminus \{l_{q}\}\qquad\text{for}\qquad  L=(l_1,\ldots,l_{|L|}).
  \end{eqnarray*}
Moreover, we  write  $\widehat I:= I\cup\{k-l\}$ for  $k-l\not\in I,$
and $\widehat J:= J\cup\{k-l\}$ for  $k-l\not\in J.$
 The coefficient of $d\zeta_{I\cup\{k-l\}} \wedge d\bar \zeta_J\wedge dw_K\wedge d\bar w_L$ in $dT$ is  given by 
 \begin{equation*}
  {\partial  T_{I,J;K,L} \over \partial \zeta_{k-l}}+\sum\limits_{q=1}^{|I|}  \pm {\partial  T_{I(q),J;K,L}\over \partial \zeta_{i_q}} 
  +\sum\limits_{q=1}^{|K|}  \pm {\partial  T_{\widehat I,J;K(q),L}\over \partial w_{k_q}}  .
 \end{equation*}
This  coefficient is zero since $dT=0.$ Moreover, we have
\begin{eqnarray*}
T_{I(q),J;K,L}(\zeta',t\zeta_{k-l},w)&=& t^{-1}T^t_{I(q),J;K,L}(\zeta,w)\quad\text{for}\quad 1\leq q\leq |I|;\\  
 T_{\widehat I,J;K(q),L}(\zeta',t\zeta_{k-l},w)&=& t^{-1}T^t_{\widehat I,J;K(q),L}(\zeta,w)\quad\text{for}\quad 1\leq q\leq |K|.
\end{eqnarray*}
Using this we infer from  the last equality that
\begin{equation*}
  {\partial  T_{I,J;K,L} \over \partial \zeta_{k-l}}={1\over t}\sum\limits_{q=1}^{|I|}  \pm {\partial  T_{I(q),J;K,L}\over \partial \zeta_{i_q}} 
  +{1\over t} \sum\limits_{q=1}^{|K|}  \pm {\partial  T_{\widehat I,J;K(q),L}\over \partial w_{k_q}}  .
 \end{equation*}
 Inserting     this  into \eqref{e:deri_f_I,J;K,L_prem} and performing an integration  by parts, it follows that
 \begin{equation}\label{e:deri_f_I,J;K,L}
 \begin{split}
 {\partial  f_{I,J;K,L} \over \partial t}&={1\over t}\sum_{q=1}^{|I|} \pm\int_U \zeta_{k-l} T^t_{I(q),J;K,L}(\zeta',t\zeta_{k-l},w) {\partial \phi(\zeta,w)\over  \partial \zeta_{i_q}}d\Leb(\zeta,w)\\
 &+{1\over t}\sum_{q=1}^{|K|} \pm\int_U \zeta_{k-l} T^t_{\widehat I,J;K(q),L}(\zeta',t\zeta_{k-l},w) {\partial \phi(\zeta,w)\over  \partial w_{k_q}}d\Leb(\zeta,w).
 \end{split}
 \end{equation}
 Similarly, we obtain  the following conjugate expression 
 \begin{equation}\label{e:deri_f_I,J;K,L-conj}
 \begin{split}
 {\partial  f_{I,J;K,L} \over \partial\bar t}&={1\over \bar t}\sum_{q=1}^{|I|} \pm\int_U \bar \zeta_{k-l} T^t_{I(q),J;K,L}(\zeta',t\zeta_{k-l},w) {\partial \phi(\zeta,w)\over  \partial \bar \zeta_{i_q}}d\Leb(\zeta,w)\\
 &+{1\over\bar  t}\sum_{q=1}^{|K|} \pm\int_U  \bar \zeta_{k-l} T^t_{ I,\widehat J;K,L(q)}(\zeta',t\zeta_{k-l},w) {\partial \phi(\zeta,w)\over  \partial \bar w_{k_q}}d\Leb(\zeta,w).
 \end{split}
 \end{equation}
 Differentiating  both sides of equality  \eqref{e:deri_f_I,J;K,L-conj}  with respect to  $t$ and using  equality \eqref{e:deri_f_I,J;K,L}, we obtain that
  \begin{equation}\label{e:Hess_f_I,J;K,L}
 \begin{split}
 {\partial^2  f_{I,J;K,L} \over \partial t\partial \bar t}&={1\over |t|^2}\sum_{1\leq q\leq |I|,\ 1\leq q'\leq |J|} \pm\int_U |\zeta_{k-l}|^2 T^t_{I(q),J(q');K,L}(\zeta',t\zeta_{k-l},w) {\partial^2 \phi(\zeta,w)\over  \partial \zeta_{i_q}    \partial \bar \zeta_{j_{q'}} }d\Leb(\zeta,w)\\
 &+{1\over |t|^2}\sum_{1\leq q\leq |K|,\ \ 1\leq q'\leq |J| } \pm\int_U |\zeta_{k-l}|^2 T^t_{\widehat I,J(q');K(q),L}(\zeta',t\zeta_{k-l},w) 
 {\partial^2 \phi(\zeta,w)\over  \partial w_{k_q}    \partial \bar \zeta_{j_{q'}}}d\Leb(\zeta,w)\\
  &+{1\over |t|^2}\sum_{1\leq q\leq |I|,\ 1\leq q'\leq |L|} \pm\int_U |\zeta_{k-l}|^2 T^t_{I(q),\widehat J;K,L(q')}(\zeta',t\zeta_{k-l},w) {\partial^2 \phi(\zeta,w)\over  \partial \zeta_{i_q}    \partial \bar w_{l_{q'}} }d\Leb(\zeta,w)\\
 &+{1\over |t|^2}\sum_{1\leq q\leq |K|,\ \ 1\leq q'\leq |L| } \pm\int_U |\zeta_{k-l}|^2 T^t_{\widehat I,\widehat J;K(q),L(q')}(\zeta',t\zeta_{k-l},w) 
 {\partial^2 \phi(\zeta,w)\over  \partial w_{k_q}    \partial \bar w_{l_{q'}}}d\Leb(\zeta,w)
 .
 \end{split}
 \end{equation}
 Note that the  function $\phi$ and  its derivatives  are  bounded on $U.$
 Since  $k-l\in I(q)$  for $1\leq q\leq |I|$  and $k-l\in J(q)$ for   $1\leq q\leq |J|,$ 
 the  first estimate of the lemma  follows  by applying Lemma \ref{L:coeff-Tt} to  equalities \eqref{e:deri_f_I,J;K,L} and \eqref{e:deri_f_I,J;K,L-conj}.
 Similarly,  the  second estimate of the lemma  follows  by applying Lemma \ref{L:coeff-Tt} to  equalities \eqref{e:Hess_f_I,J;K,L}.
 \endproof

 Recall from Blel-Demailly-Mouzali \cite{BlelDemaillyMouzali} the  following   two key lemmas.
 
 \begin{lemma}\label{L:BDM-deriv}  {\rm  \cite[Lemma 3.5]{BlelDemaillyMouzali}}
 Let $f$ be a  complex-valued  function of class $\Cc^1$ defined on  the punctured disc  $\{t\in\C: 0<|t|<\bfr\}.$
 Assume that there exists a  measurable function $u:\ (0,\bfr]\to \R^+$  such that
  \begin{equation*}
   |df(t)|\leq u(|t|)\qquad\text{and}\qquad  \int_0^\bfr  u(s)ds<\infty.
  \end{equation*}
Then  $f(t)$ admits a limit as $\C\ni t\to 0.$
 \end{lemma}

 \begin{lemma}\label{L:BDM-Hess}  {\rm  \cite[Lemma 3.6]{BlelDemaillyMouzali}}
 Let $f$ be a  complex-valued  function of class $\Cc^2$ defined on  the punctured disc  $\{t\in\C: 0<|t|<\bfr\}.$
 Assume that there exists a  measurable function $u:\ (0,\bfr]\to \R^+$  such that
  \begin{equation*}
   |\Delta f(t)|\leq u(|t|)\qquad\text{and}\qquad  \int_0^\bfr s|\log{s}| u(s)ds<\infty,
  \end{equation*}
  where $\Delta$ is  the Laplacian operator in $\C.$
Then  $f(t)$ admits a limit as $\C\ni t\to 0.$
 \end{lemma}

 \proof[End of the proof of Theorem \ref{T:main_1}]
 We only need to prove that $T_\lambda$ converge weakly  locally   in $\E|_B\setminus B$ as $|\lambda|$ tends to infinity.
 Let $U$ be  a  small relatively compact open set in  $\E|_B\setminus B$ so that we can use the  coordinate  system $y=(z,w)$ on $U.$
 Let $\phi\in\Cc^\infty_0(U)$ be a  test  function.  Let $I,J\subset\{1,\ldots,k-l\}$ and $K,L\subset \{1,\ldots,l\}.$
 We  need to show that each function
 $f_{I,J;K,L}$ defined   by \eqref{e:f_I,J;K,L}
 admits  a  limit as $\C\ni t$ tends to $0.$   Consider two cases.

 \noindent{\bf Case  $k-l\in I$ or $k-l\in J.$}
 
 By Lemma \ref{L:coeff-Tt},  the coefficient  $T^t_{I,J;K,L}$ tends in mass to $0.$
 So using  formula  \eqref{e:f_I,J;K,L}, we see that $f_{I,J;K,L}$  
also  tends   to $0$ as  $\C\ni t$ tends to $0.$

 \noindent{\bf Case  $k-l\not\in I$ and  $k-l\not\in J.$}
 
 Consider  first the subcase where   condition (a) is fulfilled.  Then we infer that $\int_0^\bfr  {\sqrt{\kappa(s)}\over s}<\infty.$
 On the other hand, by  Lemma \ref{L:f_I,J;K,L}, we have that
 \begin{equation*}
 \big| {\partial  f_{I,J;K,L} \over \partial t}\big|\leq c_1{\sqrt{\kappa(|t|)}\over |t|}.
 \end{equation*}
 So applying  Lemma \ref{L:BDM-deriv} to $f:=f_{I,J;K,L}$ and $u(s):=  {\sqrt{\kappa(s)}\over s}$ for $s\in (0,\bfr]$
 yields that  $f_{I,J;K,L}$  
 admits  a  limit as $\C\ni t$ tends to $0.$

 Consider the  last subcase where condition (b) is  fulfilled. 
 Then we infer that $\int_0^\bfr  {\kappa^\bullet(s)\over s}ds<\infty.$
 On the other hand, by  Lemma \ref{L:f_I,J;K,L}, we have that
 \begin{equation*}
 \qquad  \big| {\partial^2  f_{I,J;K,L} \over
 \partial t \partial\bar  t}\big|\leq c_2 {\kappa(|t|)\over |t|^2}
 \end{equation*}
 Observe that
 \begin{equation*}
  \int_0^\bfr  {\kappa^\bullet(s)\over s}ds=  \int_0^\bfr  {\sum_{n=0}^\infty  \kappa(2^{-n}s)\over s}ds
  =\sum_{n=0}^\infty \int_0^{2^{-n} \bfr} {\kappa(s)\over s} ds =
  \int_0^{ \bfr}\left\lbrack  {\log {(2\bfr/s)}\over \log 2} \right\rbrack {\kappa(s)\over s} ds,
 \end{equation*}
where  $[t]$ denotes the  integer part of $t\in\R.$ Hence, we  infer that
\begin{equation*}
 \int_0^{ \bfr}{  |\log{ s}| \kappa(s)\over s} ds<\infty.
 \end{equation*}
So applying  Lemma \ref{L:BDM-Hess} to $f:=f_{I,J;K,L}$ and $u(s):=  {\kappa (s)\over s^2}$ for $s\in (0,\bfr]$
 yields that  $f_{I,J;K,L}$  
 admits  a  limit as $\C\ni t$ tends to $0.$
 \endproof

 \proof[End of the proof of Theorem \ref{T:main_2}]
 Fix $i\in I.$
 Consider the functions $\kappa^{(i)},\ \kappa^\bullet:\ (0,\bfr]\to \R^+$ given by
\begin{equation}\label{e:two-kappa-i}
 \kappa^{(i)}(r):=\sum_{j=\lowm}^\upm\kappa(T,B_i\cap B,r/2,r)\quad \text{and}\quad
  \kappa^{(i),\bullet}(r):=\sum_{j=\lowm}^\upm\kappa^{(i),\bullet}(T,B_i\cap B,r)\quad \text{for}\quad  r\in (0,\bfr].
\end{equation}
If condition  (a-i)  of Theorem \ref{T:main_2} is satisfied,  then by Proposition \ref{P:main-2},
  $$\int_0^{r_0} {\sqrt{\Kc_{j,k-p-j}(T, B_i\cap B,r/2,r)}\over  r}dr <\infty\qquad \text{for all}\qquad \lowm\leq j\leq \upm;$$
and hence  we infer that
$$
\int_0^{r_0}{\sqrt{\kappa^{(i)}(r)}\over r}dr <\infty.
$$
If condition  (b-i)  of Theorem \ref{T:main_2} is satisfied,  then by Proposition \ref{P:main-2},
  $$\int_0^{r_0} {\Kc_{j,k-p-j}(T, B_i\cap B,0,r)\over  r}dr <\infty\qquad \text{for all}\qquad \lowm\leq j\leq \upm;$$
and hence  we infer that
$$
\int_0^{r_0}{\kappa^{(i),\bullet}(r)\over r}dr <\infty.
$$
Using these estimates we argue as in the proof of Theorem \ref{T:main_1}. Consequently, we can show that $T_\lambda\to  T$ on $\E_{B_i\cap B}$ as $\lambda\to\infty.$ 
Since the  $B_i$'s  form an open cover of $\overline B,$ it follows that  $T_\lambda\to  T$ on $\E_{ B}$ as $\lambda\to\infty.$ 
 \endproof

\section{ Cylindrical cone of a  complex analytic set}\label{S:Cone}

Using the  construction of the cylindrical cones in \cite{Nguyen24b},  we  prove the last main result (Theorem \ref{T:main_3}) in this  section.

Let $(X,\omega)$ be a  K\"ahler manifold of dimension $k,$ and let  $V\subset X$ be a submanifold of dimension $l, $
and let $B\subset V$ be a relatively compact piecewise $\Cc^2$-smooth open subset. Let $\pi:\ \E\to V$ be the  normal vector  bundle to $V$ in $X.$  
 Let $S$ be  a  complex analytic set of pure codimension $p$ in $X$ such that $S\cap B\Subset B.$ 
 First  we  recall the construction of the  {\it cylindrical cone} $\Cf(S)$ of $S$   described  in \cite{Nguyen24b}, this is a $V$-conic analytic set on 
 $\E|_B.$ We consider two cases.

 \noindent  {\bf  Local  case:}

 In this case 
  $(S,0)$ is a germ of complex analytic set of pure codimension $p$ in $\C^k.$
and  $(V,0)$ is   a germ  of a submanifold of dimension $l$  in $\C^k.$  We do not need to assume that  $S\cap B\Subset B.$
This  case is very intuitive and  it helps us to understand the geometric meaning of the cylindrical cone.
Let $y=(z,w)\in\C^{k-l}\times \C^l=\C^k$ be the coordinates of $\C^k.$ 
We may assume without loss of generality that $$V=\{0\}\times \C^l=\{(0,w):\ w\in\C^l\}.$$ 
We may also  suppose  that $S$ is  defined  by  the  equations $f_1(z,w)=\ldots=f_N(z,w)=0$
on an open neighborhood $\Omega$ of $0$ in $\C^k,$ where $f_1,\ldots,f_N$ are holomorphic  functions on $\Omega.$
  Consider  the open set
 $$
 U:=\left\lbrace (t,z,w)\in\C\times\C^{k-l}\times \C^l:\  (tz,w)\in\Omega \right\rbrace\subset \C^{k+1}.
 $$
Consider also  the analytic set $E$ in $U$ defined by  the  equations $f_1(tz,w)=\ldots=f_N(tz,w)=0.$
It is 
 the union of $\{0\}\times\C^{l}\times \C^{k-l}$ and the dilation $\{t\}\times t^{-1}S$ for $t\in\C^*.$   Let $E^*$ be the union of irreducible components
 of $E$ which are not contained in $\{0\}\times\C^k.$ Every connected component $E_j$ of $E$ is of dimension $k-p+1.$  Moreover,
 $E_j\setminus (\{0\}\times\C^k)$ is dense in $E_j.$  Hence, $E^\star=\bigcup E_j$ is  the closure in $\C^{k+1}$ of
 $$
 E\setminus (\{0\}\times\C^k)=\bigcup_{\lambda \in\C^*}\{1/\lambda\}\times \lambda S.
 $$
 The  {\it  cylindrical cone}  $\Cf(S)$ is  the analytic set of dimension $k-p$ in $\C^k$ defined by
 \begin{equation}\label{e:CC}
 \{0\}\times \Cf(S)=E^*\cap (\{0\}\times \C^k)\subset \C^{k+1}.
\end{equation} 
 As observed  in  \cite{Nguyen24b},    $\Cf(S)$ is exactly the set of limits of the sequences  $(\lambda_n z^{(n)},w^{(n)}),$
 where $y^{(n)}= (z^{(n)},w^{(n)})\in S$ and  $\lambda_n\in\C^*$ tends to $\infty$ as $n$ tends to infinity.

 \noindent  {\bf  Global  case:}
 
 This is   the general  case and  we   need to assume that  $S\cap B\Subset B.$
 Let $\tau$ be an admissible map along $B.$
 The  {\it  cylindrical cone}  $\Cf(S)$ is  the analytic set of dimension $k-p$ in $\E|_B $ defined by
 \begin{equation}\label{e:CC-bis}
 \{0\}\times \Cf(S)=E^*\cap (\{0\}\times \E|_B).
\end{equation} 
 Observe that   $\Cf(S)$ is exactly the set of limits of the sequences  $(\lambda_n z^{(n)},w^{(n)}),$
 where $y^{(n)}= (z^{(n)},w^{(n)})\in \tau(S)$ and  $\lambda_n\in\C^*$ tends to $\infty$ as $n$ tends to infinity.

 Let $[S]$ be the current of integration on $S.$ 
 We infer  from the above  paragraph and  the asssumption  $S\cap B\Subset B$ that 
 \begin{equation}\label{e:Cf(S)_ver-empty} \Cf(S)\cap  \partial_\ver\Tube(B,1).
 \end{equation}
   Consider the function $S_\reg\ni (z,w)\mapsto \varphi(z,w).$
By Sard's theorem, the set of all  critical values $r$ of  this  function, denoted by $\Ec,$  is  at most countable.
 Since for $\lambda\in\C^*,$ $(A_\lambda)_*[\tau_*S]=[A_\lambda(\tau(S))],$  we have
 for $0<r<\bfr$ with $r\not\in\Ec$  and $\lowm\leq j\leq \upm$ that 
 \begin{equation*}
\nu_j([S],B,r,\tau,h)=\int_{(A_{r^{-1}})(\tau(S))\cap \Tube(B,1)}\pi^*\omega^j\wedge \beta^{k-p-j}.
 \end{equation*}
  By  Stokes' theorem applied to all non-critical values $r$ of  the function $S_\reg\ni (z,w)\mapsto \varphi(z,w),$ and by using
  the  expression $\beta=\ddc\varphi$  given by \eqref{e:alpha-beta-spec}, 
the RHS is   equal to 
\begin{equation*}
   \int_{A_{r^{-1}}(\tau(S))\cap \partial \Tube(B,1)}\pi^*\omega^j\wedge
  \beta^{k-p-j-1}\wedge \dbar\varphi.
  \end{equation*}
Therefore, we infer that
\begin{equation}\label{e:nu_j-S}
\nu_j([S],B,r,\tau,h)=\big(\int_{A_{r^{-1}}(\tau(S))\cap \partial_\hor \Tube(B,1)}+  \int_{A_{r^{-1}}(\tau(S))\cap \partial_\ver\Tube(B,1)}\big)\pi^*\omega^j\wedge
  \beta^{k-p-j-1}\wedge \dbar\varphi.
 \end{equation} 
 As $\R^+\subset \C$ and  $\Tube(B,1)\subset\E|_B,$ consider the 
following  subset of $\C\times \E|_B:$ 
 \begin{equation*}
  M:=E^{*}\cap (\R^+\times \partial\Tube(B,1)):=M_\hor\cup M_\ver,
 \end{equation*}
 where
  \begin{equation*}
  M_\hor:=E^{*}\cap (\R^+\times \partial_\hor\Tube(B,1))\quad\text{and}\quad   M_\ver:=E^{*}\cap (\R^+\times \partial_\ver\Tube(B,1)).
 \end{equation*}
For $0\leq r_1<r_2\leq \bfr,$  as  $(r_1,r_2)\subset \C,$ consider following  subset of $\C\times \E|_B:$ 
\begin{equation*}
  M(r_1,r_2)=E^{*}\cap ( (r_1,r_2)\times \partial\Tube(B,1))=M\cap ( (r_1,r_2)\times\E|_B ), 
 \end{equation*}
 and  define in the same way the two sets  $M_\hor(r_1,r_2),$ $ M_\ver(r_1,r_2)$ using   $M_\hor, $  $M_\ver$ instead of $M.$
We have  for $r\not\in\Ec,$  $$\dim_\R A_{r^{-1}} (\tau(S))\cap   \partial_\hor\Tube(B,1)=2k-2p-1 . $$ Moreover, since  $S\cap B\Subset B,$  we may  fix  $\bfr>0$ so small enough   that  for $0<r\leq\bfr,$
\begin{equation}\label{e:A_r_-1(tau(S))empty} A_{r^{-1}}( \tau(S))\cap   \partial_\ver\Tube(B,1)=\varnothing.
\end{equation}
Note  that
\begin{eqnarray*}
M_\hor(r_1,r_2)&=&\{0\}\times (\Cf(S)\cap  \partial_\hor\Tube(B,1))\cup \bigcup_{r\in(r_1,r_2)} \{r\}\times  ( A_{r^{-1}}( \tau(S))\cap  \partial_\hor\Tube(B,1)),\\
M_\ver(r_1,r_2)&=&\{0\}\times (\Cf(S)\cap  \partial_\ver\Tube(B,1))\cup \bigcup_{r\in(r_1,r_2)} \{r\}\times  ( A_{r^{-1}}(\tau( S))\cap  \partial_\ver\Tube(B,1)),\\
&=& \bigcup_{r\in(r_1,r_2)} \{r\}\times  ( A_{r^{-1}}(\tau( S))\cap  \partial_\ver\Tube(B,1)),
\end{eqnarray*}
where the last equality holds by \eqref{e:Cf(S)_ver-empty}.
We infer from  the first equality of the last line  that $M_\hor$ is  of pure real dimension $2k-2p.$ For $r_1,r_2\not\in \Ec,$   the horizontal boundary $\partial_\hor M(r_1,r_2)$
is identified   to 
$$
 \bigcup_{j=1,2} \{r_j\}\times  ( A_{r_j^{-1}} (\tau(S))\cap  \partial_\hor\Tube(B,1)),
$$
and  it is  smooth at every regular point of  $A_{r_j^{-1}}(\tau( S)).$ 
We also   infer from    the last equality  
for $M_\ver(r_1,r_2)$ and 
from
\eqref{e:A_r_-1(tau(S))empty} that
\begin{equation}\label{e:M_ver-empty}
M_\ver(r_1,r_2)=\varnothing \qquad\text{for}\qquad 0\leq r_1<r_2\leq \bfr.
\end{equation}
Formula  \eqref{e:nu_j-S} and Stokes' theorem  give that
\begin{equation*}
  \nu_j([S],B,r_2,\tau,h)  -\nu_j([S],B,r_1,\tau,h)=\int_{M(r_1,r_2)}\pi^*\omega^j\wedge \beta^{k-p-j}.
\end{equation*}
This  formula  remains valid by continuity  for all values $0\leq r_1<r_2.$   Wirtinger's inequlity shows that $(c_1\pi^*\omega+c_2\beta)^{k-p}$
is bounded  from above by a constant times the Riemannian $2(k-p)$-volume of $M.$ We obtain in particular for $r_1=0$ and $r_2=r\in (0,\bfr]$: 
\begin{equation}\label{e:nu_j-bounded-by-vol}
  |\nu_j([S],B,r,\tau,h)  -\nu_j([S],B,h)|\leq c  \cdot \volume_{2(k-p)} (M(0,r))\qquad\text{for}\qquad  \lowm\leq j\leq \upm.
\end{equation}
The  following result is  needed.
\begin{lemma}\label{L:Volume} {\rm (Blel-Demailly-Mouzali \cite[Lemme 5.5]{BlelDemaillyMouzali})}
 Let $A$ be a real analytic set  of dimension $\leq n$ in an open set $\Omega\subset \R^N$ and let  $g$ be a  real analytic function  on $\Omega.$ Consider the set
 $$
 A(r):=\left\lbrace  y\in A:\   0< g(y)<r \right\rbrace.
 $$
 Then  for every compact set $K\subset \Omega,$ there are constants $C,\rho>0$  such that
 $$
 \volume_{n} (A(r)\cap K)\leq cr^\rho.
 $$
\end{lemma}

 \proof[End of the proof of Theorem \ref{T:main_3}]
      Consider  the closure $A$  of $M(0,\bfr).$   Since we know  by  \eqref{e:M_ver-empty} that
$M_\ver(0,\bfr)=\varnothing,$  $A$ is   real analytic.  For $r\in (0,\bfr],$  The set $A(r)$ of Lemma \ref{L:Volume}  satisfies
$A(r)=M_\hor(0,r)$ and  by that lemma, there is  a constant $\rho>0$ such that   $\volume_{2k-2p} (A(r))\leq cr^\rho.$ 
 This, combined with estimate \eqref{e:nu_j-bounded-by-vol}, implies  the  result.
 \endproof

 \begin{remark}\rm
  The reader   may consult  \cite{Nguyen24b} for more  extensive  results of the generalized Lelong numbers  for the  current of integration
  over an analytic  set. The treatment  therein   emphasizes     the link between    algebraic geometry  and complex   geometry.
 \end{remark}

\end{document}